\documentclass{amsart}
\usepackage{amssymb}

\newtheorem{prop}{Proposition}
\newtheorem{thm}[prop]{Theorem}
\newtheorem{cor}[prop]{Corollary}
\newtheorem{lem}[prop]{Lemma}
\newtheorem{conj}[prop]{Conjecture}
\theoremstyle{definition}
\newtheorem{ex}[prop]{Example}
\newtheorem{rem}[prop]{Remark}

\begin{document}

\title[]{Graded characters of modules \\
supported in the closure of \\
a nilpotent conjugacy class}
	
\author[M. Shimozono]{Mark Shimozono}
\address{Dept.of Mathematics \\ Virginia Tech \\ Blacksburg, VA }
\email{mshimo@math.vt.edu}
\thanks{First author partially supported by a postdoctoral supplement
under Richard Stanley's NSF Research Grant DMS-9500714.}

\author[J. Weyman]{Jerzy Weyman}
\address{Dept.\ of Mathematics \\ Northeastern University \\
Boston, MA}
\email{weyman@neu.edu}
\thanks{Second author partially supported by NSF Research Grant
DMS-9403703.}

\begin{abstract} This is a combinatorial study of
the Poincar\'e polynomials of isotypic
components of a natural family of graded $GL(n)$-modules
supported in the closure of a nilpotent conjugacy class.
These polynomials generalize the Kostka-Foulkes polynomials
and are $q$-analogues of Littlewood-Richardson coefficients.
The coefficients of two-column Macdonald-Kostka polynomials
also occur as a special case.
It is conjectured that these $q$-analogues  are the generating function
of so-called catabolizable tableaux with the charge statistic of
Lascoux and Sch\"utzenberger.  A general approach for a proof is
given, and is completed in certain special cases including the
Kostka-Foulkes case.  Catabolizable tableaux are used to prove
a characterization of Lascoux and Sch\"utzenberger for the
image of the tableaux of a given content under the standardization
map that preserves the cyclage poset.
\end{abstract}

\maketitle

\newcommand{\la}{\lambda}
\newcommand{\dom}{\trianglerighteq}
\newcommand{\C}{\mathbb{C}}
\newcommand{\N}{\mathbb{N}}
\newcommand{\Z}{\mathbb{Z}}
\newcommand{\LR}{\mathrm{LR}}
\newcommand{\rp}{{r\!+\!1}}
\newcommand{\ul}{\underline}
\newcommand{\tK}{\widetilde{K}}
\newcommand{\K}{K}
\newcommand{\Kn}{\sim_K}

\newcommand{\muhat}{{\widehat \mu}}
\newcommand{\Sh}{\widehat{S}}
\newcommand{\Th}{\widehat{T}}
\newcommand{\Uh}{\widehat{U}}
\newcommand{\gammahat}{{\widehat \gamma}}
\newcommand{\etahat}{{\widehat \eta}}
\newcommand{\Rhat}{{\widehat R}}
\newcommand{\shape}{\mathrm{shape}}
\newcommand{\charge}{\mathrm{charge}}
\newcommand{\cocharge}{\mathrm{cocharge}}
\newcommand{\sign}{\mathrm{sign}}
\newcommand{\weight}{\mathrm{weight}}
\newcommand{\cattype}{\mathrm{cattype}}
\newcommand{\Lie}{\mathrm{Lie}}
\newcommand{\ch}{\mathrm{ch}}
\newcommand{\GL}{\mathrm{GL}}
\newcommand{\gl}{\mathrm{gl}}
\newcommand{\LLL}{\mathcal{L}}
\newcommand{\gggg}{\mathfrak{g}}
\newcommand{\hhh}{\mathfrak{h}}
\newcommand{\ppp}{\mathfrak{p}}
\newcommand{\bbb}{\mathfrak{b}}
\newcommand{\MM}{\mathcal{M}}
\newcommand{\OOO}{\mathcal{O}}
\newcommand{\RR}{\mathcal{R}}
\newcommand{\SSS}{\mathcal{S}}
\newcommand{\TTT}{\mathcal{T}}
\newcommand{\diag}{\mathrm{diag}}
\newcommand{\Hom}{\mathrm{Hom}}
\newcommand{\trace}{\mathrm{tr}}
\newcommand{\ev}{\mathrm{ev}}
\newcommand{\cat}{\mathrm{cat}}
\newcommand{\Cat}{\mathrm{Cat}}
\newcommand{\CCat}{\mathrm{CCat}}
\newcommand{\sq}{\times}
\newcommand\Roots{\mathrm{Roots}}
\newcommand{\inner}[2]{\langle #1,#2\rangle}
\newcommand{\Fdot}{F_\cdot}

\numberwithin{equation}{section}

\section{Introduction}

For a partition $\mu$ of $n$ let $X_\mu$ be the Zariski closure of the
conjugacy class in $gl(n,\C)$ of the nilpotent Jordan matrix with blocks
of sizes given by the conjugate or transpose partition $\mu^t$ of $\mu$.
Since $X_\mu$ is a cone that is stable under the action of $GL(n,\C)$ given
by matrix conjugation, its coordinate ring $\C[X_\mu]$ is graded and affords
an action of $GL(n,\C)$ that respects the grading.  The natural category
of modules over $\C[X_\mu]$ is the family of finitely generated graded
$\C[X_\mu]$-modules (that is, $\C[gl(n)]$-modules supported in $X_\mu$)
that afford the action of $GL(n,\C)$ that is compatible
with the graded $\C[X_\mu]$-module structure.

For any permutation $\eta$ of the parts of $\mu$, there is a Springer
desingularization $Z(\eta)\rightarrow X_\mu$.  Corresponding to any
$GL(n)$-weight $\gamma\in\Z^n$, there is an $\OOO_{Z(\eta)}$-module
$\MM_{\eta,\gamma}$ with Euler characteristic $\chi_{\eta,\gamma}$,
which can be viewed as an element of the Grothendieck group
$K_0'(\C[X_\mu])$ of the aforementioned category of $\C[X_\mu]$-modules.
In \cite{KKSW} it is shown that the classes of the $\chi_{\eta,\gamma}$
graded $GL(n)$-modules $\chi_{\eta,\gamma}$ generate the group
$K_0'(\C[X_\mu])$.

The purpose of this paper is to investigate the combinatorial
properties of the family of polynomials $\K_{\la,\gamma,\eta}(q)$
given by the isotypic components of the virtual graded
$GL(n)$-modules $\chi_{\eta,\gamma}$.

The polynomials $\K_{\la,\gamma,\eta}(q)$ are
$q$-analogues of Littlewood-Richardson (LR) coefficients.
Special cases include the Kostka-Foulkes polynomials
(Lusztig's $q$-analogues of weight multiplicities in type A)
and coefficients of two-column Macdonald-Kostka polynomials.

Many formulas involving the Kostka-Foulkes polynomials
have suitable generalizations for the polynomials $\K_{\la,\gamma,\eta}(q)$,
such as the $q$-Kostant partition function formula \cite{Mac}
and Morris' recurrence \cite{Mo}.  We derive these formulas directly
from the definition of the twisted modules $\MM_{\eta,\gamma}$.

The Kostka-Foulkes polynomials also have several combinatorial
descriptions, including two beautiful formulas of Lascoux and
Sch\"utzenberger involving the $q$-enumeration of two sets of 
tableaux \cite{La} \cite{LS1}.  The main focus of this paper
is a common generalization of these formulas.
We define the notion of a \textit{catabolizable tableau}
and give a conjectural interpretation of $\K_{\la,\gamma,\eta}(q)$ as the
generating function over such tableaux with the \textit{charge}
statistic (Conjecture \ref{cat conj}).
A general approach for a proof is given, using a sign-reversing
involution that cancels terms from the generalized Morris recurrence.
The missing ingredient in the general case is to
prove that the involution preserves catabolizability in a certain sense.
In special cases this can be shown, so the conjecture holds in those
cases.

We also develop properties of catabolizable tableaux and their
intimate relationship with the cyclage poset \cite{La} \cite{LS1}.
As an application, we supply a proof of a formula of Lascoux \cite{La}
for the cocharge Kostka-Foulkes polynomials.

In the special case corresponding to LR coefficients associated
with products of rectangular shapes, the polynomials $K_{\la,\gamma,\eta}(q)$
seem to coincide with yet another $q$-analogue of LR coefficients
given by the combinatorial objects known as rigged configurations.
This connection is pursued in \cite{KS}.

Lascoux, Leclerc and Thibon have defined a $q$-analogue of certain
LR coefficients using ribbon tableaux \cite{LLT}.  This family
of polynomials appears to contain the polynomials
$K_{\la,\gamma,\eta}$ as a subfamily, but the reason for this is
unclear.

\section{The polynomials $\K_{\la,\gamma,\eta}(q)$}

The first goal is to derive an explicit formula for the
polynomials $\K_{\la,\gamma,\eta}(q)$.  By definition this polynomial
is the coefficient of the irreducible character $s_\la(x)$ of highest
weight $\la$ in the formal character $H_{\eta,\gamma}(x;q)$
of the virtual graded $GL(n)$-module $\chi_{\eta,\gamma}$.  
A nice bialternant formula for $H_{\eta,\gamma}(x;q)$ is obtained
by expressing the Euler characteristic $\chi_{\eta,\gamma}$ in terms
of Euler characteristics of $GL(n)$-equivariant line bundles over the flag
variety and applying Bott's theorem for calculating the latter.
Some readers may prefer to take the
formulas \eqref{gen func def} and \eqref{poly def} as the definitions of
$H_{\eta,\gamma}(x;q)$ and $K_{\la,\eta,\gamma}(x;q)$ respectively.

The rest of this section derives various properties
of the polynomials $K_{\la,\eta,\gamma}(x;q)$, including positivity
conditions, specializations to known families of polynomials,
the interpretation as a $q$-analogue of an LR coefficient,
and various symmetry and monotonicity properties.

\subsection{The Poincar\'e polynomial $K_{\la,\gamma,\eta}(q)$}
\label{Euler def}
We review the definition of $\chi_{\eta,\gamma}$ given in
\cite{KKSW} and derive a formula for
its graded character $H_{\eta,\gamma}(x;q)$ and its coefficient
polynomials $K_{\la,\gamma,\eta}(x;q)$.

Let $\mu$ be a partition of a fixed positive integer $n$
and $X_\mu$ the nilpotent adjoint orbit closure defined in the introduction.
Let $\eta=(\eta_1,\dots,\eta_t)$ be a reordering of the parts of the
partition $\mu$.  For each such $\eta$ 
there is a Springer desingularization $q_\eta:Z_\eta\rightarrow X_\mu$
defined as follows.  Consider the variety of partial flags of dimensions
given by the sequence $0 = d_t < d_{t-1} <\dots <d_1<d_0=n$, where
\begin{equation*}
   d_i = n - (\eta_1+\dots+\eta_i)\qquad\text{ for $0\le i\le t$,}
\end{equation*}
whose typical element is 
\begin{equation*}
  \Fdot = (0=F_{d_t}\subset F_{d_{t-1}}\subset\dots\subset F_{d_0}=\C^n)
\end{equation*}
where $F_j$ is a subspace of dimension $j$.  This flag variety is
realized by the homogeneous space $G/P$, where $G=GL(n,\C)$ and
$P=P_\eta$ is the parabolic subgroup that stabilizes the partial flag
whose subspaces have the form $\C^{d_i}$ for $0\le i\le t$, where $\C^j$
denotes the span of the \textit{last} $j$ standard basis vectors of the
standard left $G$-module $\C^n$, viewed as the space of column vectors.
In other words, $P$ is the subgroup of lower block triangular matrices
with block sizes given in order by $\eta_1$ through $\eta_t$.
Let $\gggg=Lie(G)$ and $\ppp=\ppp_\eta=Lie(P)$.

Define the incidence variety
\begin{equation*}
  Z = Z_\eta := \{(A,\Fdot)\in \gggg\times G/P \mid A F_{d_{i-1}} \subset
	F_{d_i} \text{for $1\le i\le t$}\}.
\end{equation*}
Let $q:Z\rightarrow \gggg$ and $p:Z\rightarrow G/P$ be the restriction
to $Z$ of the first and second projections of $\gggg\times G/P$.
The map $q$ is a desingularization of its image
$X_\mu$ \cite{He}.

Next we recall the definition of the $\OOO_Z$-modules $\MM_{\eta,\gamma}$.
Let $B\subset G$ be the standard \textit{lower} triangular
Borel subgroup, $H\subset B$ the subgroup of diagonal matrices, and
$\Delta^+$ the set of positive roots, which are chosen to be in the
\textit{opposite} Borel to $B$.  Let $\bbb=\Lie(B)$ and
$\hhh=\Lie(H)$.  Let $W$ be the Weyl group, the symmetric group on
$n$ letters, which shall be identified with the permutation matrices in $G$.
The character group of $H$ (and the integral weights) may be
identified with $\Z^n$ such that $(\gamma_1,\dots,\gamma_n)\in\Z^n$ is
identified with the character
$(\diag(x_1,\dots,x_n)\mapsto \prod_{i=1}^n x_i^{\gamma_i}$
where $\diag(x_1,\dots,x_n)$ is the diagonal matrix with
diagonal entries $x_i$.
According to these conventions, a weight $\gamma\in\Z^n$ is dominant
if $\gamma_i\ge \gamma_{i+1}$ for $1\le i\le n-1$.

Let $\C_\gamma$ be the one-dimensional $B$-module of weight $\gamma$
and $\LLL_\gamma := G \times^B \C_\gamma$ the $G$-equivariant line
bundle over $G/B$ given by the orbit space
\begin{equation*}
\begin{split}
  G \times^B \C_\gamma &:= (G \times \C_\gamma)/B \\
  (g,v)b = (gb,\lambda(b^{-1}) v)
\end{split}
\end{equation*}
with bundle map $(g,v)B\mapsto gB$.
Let $\phi=\phi_\eta:G/B\rightarrow G/P$ be the canonical projection.
Define the $\OOO_Z$-module
\begin{equation*}
  \MM_{\eta,\gamma} = \OOO_Z \otimes p^* \phi_*(\LLL_\gamma).
\end{equation*}
Define the elements $\chi_{\eta,\gamma}$
and $M_{\eta,\gamma}$ in $K_0'(\C[X_\mu])$ by 
\begin{equation*}
\begin{split}
  \chi_{\eta,\gamma} &= \sum_{i\ge 0} (-1)^i
	[ \RR^i q_*(\MM_{\eta,\gamma})] \\
  M_{\eta,\gamma} &= q_*(\MM_{\eta,\gamma}).
\end{split}
\end{equation*}

Next, the element $\chi_{\eta,\gamma}$ is expressed in terms of Euler
characteristics of $G$-equivariant line bundles over $G/B$.
Since $q:Z\rightarrow \gggg$ is a morphism to an affine variety,
\begin{equation*}
	\RR^i q_* (\MM_{\eta,\gamma}) \cong
	H^i(Z,\MM_{\eta,\gamma})\qquad\text{for $i\ge0$.}
\end{equation*}
Since
\begin{equation*}
  \RR^i p_*(\MM_{\eta,\gamma}) = 0 \qquad\text{for $i>0$},
\end{equation*}
we have
\begin{equation*}
  H^i(Z,\MM_{\eta,\gamma}) \cong H^i(G/P,p_*(\MM_{\eta,\gamma}))
  \qquad\text{for $i\ge0$.}
\end{equation*}
By the projection formula,
\begin{equation*}
  p_*(\MM_{\eta,\gamma}) =
  p_*(\OOO_Z \otimes p^* \phi_* \LLL_\gamma) =
  p_*(\OOO_Z) \otimes \phi_* \LLL_\gamma.
\end{equation*}
As sheaves over $G/P$, $p_*(\OOO_Z)$ may be identified
$\TTT := G \times^P S(\gggg/\ppp)$ where $S$ is the symmetric algebra.
Thus we have
\begin{equation*}
  p_*(\OOO_Z) \otimes \phi_*\LLL_\gamma =
  \TTT \otimes \phi_*\LLL_\gamma.
\end{equation*}
Let $\TTT' := G \times^B S(\gggg/\ppp)$.  Then $\TTT'=\phi^*\TTT$ and
\begin{equation*}
  R^i \phi_*(\TTT' \otimes \LLL_\gamma) = 0 \qquad\text{for $i>0$.}
\end{equation*}
Then 
\begin{equation*}
\begin{split}
     H^i(G/B,\TTT' \otimes \LLL_\gamma)
  &= H^i(G/B,\phi^*(\TTT) \otimes \LLL_\gamma) \\
  &= H^i(G/P,\phi_*(\phi^*(\TTT) \otimes \LLL_\gamma) \\
  &= H^i(G/P,\TTT \otimes \phi_*\LLL_\gamma))
\end{split}
\end{equation*}
Putting this all together, we have
\begin{equation} \label{chi flag bundle}
  \chi_{\eta,\gamma} =
  \sum_{i\ge0} (-1)^i [ H^i(G/B,\TTT' \otimes \LLL_\gamma) ].
\end{equation}
The homogeneous degree is given by the degree in the polynomial ring
$S(\gggg/\ppp)$.  Each homogeneous component of $\TTT' \otimes \LLL_\gamma$
has a filtration whose successive quotients are $G$-equivariant line bundles
over $G/B$.  By the additivity of Euler characteristic the calculation
reduces to Bott's formula for the Euler characteristic
of a line bundle over $G/B$.

This allows the explicit calculation of the formal character
$H_{\eta,\gamma}(x;q)$ of $\chi_{\eta,\gamma}$.
The \textit{formal character} of a finite dimensional rational
$H$-module $M$ is the Laurent polynomial given by
\begin{equation*}
  \ch(M) = \trace (x|M),
\end{equation*}
the trace of the action of $x=\diag(x_1,\dots,x_n)$ on $M$.

Let $\rho=(n-1,n-2,\dots,1,0)$ and $J$ and $\pi$ the operators
on $\C[x_1,\dots,x_n][\det(x)^{-1}]$ given by
\begin{equation*}
\begin{split}
  J(f) &= \sum_{w\in W} (-1)^w w f \\
  \pi(f) &= J(x^\rho)^{-1} J(x^\rho f).
\end{split}
\end{equation*}
The operator $\pi$ is the Demazure operator
$\pi_{w_0}$ in the notation of \cite{Mac}, where $w_0$ is the
longest element in $W$.  Bott's formula states that
\begin{equation*}
  \ch(\chi(G/B,\LLL_\alpha)) = \pi(x^\alpha)
\end{equation*}
for $\alpha\in\Z^n$.  In particular, if $\alpha$ is a dominant integral
weight (resp. partition) then $\pi(x^\alpha)=s_\alpha(x)$ is the
irreducible character of highest weight $\alpha$ (resp. Schur polynomial).

Consider the weights of the adjoint action of $P$ on $\gggg/\ppp$,
which is indexed by the set
$\Roots_\eta$ of matrix positions above the
block diagonal given by the parts of $\eta$, that is,
\begin{equation*}
  \Roots_\eta = \{(i,j) \mid
	1 \le i \le \eta_1+\dots+\eta_r < j \le n
	\text{ for some $r$}\}.
\end{equation*}
\begin{ex}
\begin{enumerate}
\item If $\eta=(n)$ then $\Roots_\eta$ is empty.
\item If $\eta=(1^n)$ then $\Roots_\eta=\{(i,j):1\le i<j\le n\}$.
\end{enumerate}
\end{ex}

Keeping track of the degree in $S(\gggg/\ppp)$ by powers of the
variable $q$, let the formal power series $B_\eta(x;q)$ 
(resp. $H_{\eta,\gamma}(x;q)$) be the formal character of the graded
$B$-module $S(\gggg/\ppp)$ (resp. the $G$-module
$\chi_{\eta,\gamma}$).  By \eqref{chi flag bundle} and Bott's formula,
these can be written
\begin{equation*}
\begin{align} \label{gen func def}
  B_\eta(x;q) &= \prod_{(i,j)\in\Roots_\eta} \dfrac{1}{1-q\, x_i/x_j} \\
  H_{\eta,\gamma}(x;q) &= \pi ( x^\gamma B_\eta(x;q) ).
\end{align}
\end{equation*}
Then by definition,
\begin{equation} \label{poly def}
  H_{\eta,\gamma}(x;q) = \sum_\la \K_{\la,\gamma,\eta}(q) s_\la(x),
\end{equation}
where $\la$ runs over the dominant integral weights in $\Z^n$.
$B_{\eta}(x;q)$ and $H_{\eta,\gamma}(x;q)$
should be viewed as formal power series in $q$ with coefficients
in the ring of formal Laurent polynomials in the $x_i$.
It is shown later in Proposition \ref{q Kostant} that
$\K_{\la,\gamma,\eta}(q)$ is a polynomial with integer coefficients.

\begin{ex} Let $n=2$, $\eta=(1,1)$ and $\gamma=(0,0)$.  Then
\begin{equation*}
\begin{split}
  H_{\eta,\gamma}(x;q) &= \sum_{k\ge0} q^k s_{(k,-k)}(x_1,x_2) \\
  &= \sum_{k\ge0} q^k (x_1 x_2)^{-k} s_{(2k,0)}(x_1,x_2)
\end{split}
\end{equation*}
So $\K_{(k,-k),(0,0),(1,1)}(q) = q^k$ for all $k\in\N$.
\end{ex}

\subsection{Normalization}

Many of the polynomials $\K_{\la,\gamma,\eta}(q)$
coincide for different sets of indices.  We show that certain
simplifying assumptions can be made on the indices, and define
another notation that explicitly indicates
the connection with LR coefficients.

First, observe that $\K_{\la,\gamma,\eta}(q) = 0$
unless $|\la|=|\gamma|$, where $|\gamma| := \sum_{i=1}^n \gamma_i$.
To see this, consider Bott's formula for the operator $\pi$
acting on a Laurent monomial.  For $\alpha\in\Z^n$,
Let $\alpha^+$ be the unique dominant weight in the $W$-orbit of $\alpha$,
and let $w_\alpha\in W$ be the shortest element
such that $w \alpha^+ = \alpha$.  Then
\begin{equation}
\label{Bott formula}
  \pi(x^\alpha) =
    \begin{cases}
      0 & \text{if $\alpha+\rho$ has a repeated part} \\
      (-1)^{w_{\alpha+\rho}} s_{(\alpha+\rho)^+ - \rho}(x) & \text{otherwise}
    \end{cases}
\end{equation}
Every Laurent monomial $x^\alpha$ in $x^\gamma B_\eta(x;q)$
satisfies $|\alpha|=|\gamma|$, and if $\pi(x^\alpha)$ is nonzero
then it equals $\pm s_\la(x)$ where $|\alpha|=|\la|$, proving the
assertion.

Next, it may be assumed that the weights $\la$ and $\gamma$
have nonnegative parts.  To see this, let $\gamma+k$ denote the weight
obtained by adding the integer $k$ to each part of $\gamma$.  Since
$(x_1\dots x_n)^k$ is $W$-symmetric, it follows that
\begin{equation*} \label{shift all}
\begin{split}
  s_{\la+k}(x) &= (x_1 x_2 \dots x_n)^k s_\la(x) \\
  H_{\eta,\gamma+k}(x;q) &= (x_1 x_2 \dots x_n)^k H_{\eta,\gamma}(x;q) \\
  \K_{\la+k,\gamma+k,\eta}(q) &= \K_{\la,\gamma,\eta}(q)
\end{split}
\end{equation*}

Given the pair $(\eta,\gamma)$, let
$R=R(\eta,\gamma)=(R_1,R_2,\dots,R_t)$ be the
sequence of weights where $R_1$ is the $GL(\eta_1)$-weight
given by the first $\eta_1$ parts of $\gamma$,
$R_2$ the $GL(\eta_2)$-weight given by the next $\eta_2$ parts of $\gamma$,
etc.  Then it may be assumed that $R_i$ is a dominant weight of $GL(\eta_i)$
for each $i$.  For this, we require a few properties of the
isobaric divided difference operators, whose proofs are easily
derived from \cite[Chapter II]{Mac2}.

Let $I$ be a subinterval of the set $[n]=\{1,2,\dots,n\}$
and $\pi_I$ the isobaric divided difference operator indexed by
the longest element in the symmetric group on the set $I$.  Then
\begin{enumerate}
\item $\pi_I f = f$ provided that $f$ is symmetric in the
variables $\{x_j: j\in I\}$.
\item $\pi \pi_I = \pi$.
\end{enumerate}

Let $A_1$ be the first $\eta_1$ numbers in $[n]$, $A_2$ the
next $\eta_2$ numbers, etc.  That is,
\begin{equation}
A_i = [\eta_1+\dots+\eta_{i-1}+1,\eta_1+\eta_2+\dots+\eta_i]
\end{equation}
for $1\le i\le t$.  Using the appropriate symmetry properties of
$B_\eta(x;q)$, we have
\begin{equation*}
\begin{split}
  H_{\eta,\gamma}(x;q) &= \pi \,(x^\gamma B_\eta(x;q)) \\
  &= \pi \,\pi_{A_i} \,(x^\gamma B_\eta(x;q)) \\
  &= \pi\, B_\eta(x;q)\, \pi_{A_i}\, x^\gamma.
\end{split}  
\end{equation*}
By two applications of \eqref{Bott formula}
applied to this subset of variables, it follows that
$H_{\eta,\gamma}(x;q)$ is either zero, or (up to sign)
equal to $H_{\eta,\gamma'}(x;q)$ for some weight $\gamma'$ such that
the associated $GL(\eta_i)$-weights $R_i'$ are dominant.

\begin{rem} \label{normal parameters}
To summarize, the polynomial $\K_{\la,\gamma,\eta}(q)$
is either zero or (up to sign) equal to another such polynomial where
\begin{enumerate}
\item $\la=(\la_1,\dots,\la_n)$ is a partition with $n$ parts
(some of which may be zero);
\item The weight $R_i$ is a partition with $\eta_i$ parts
(some of which may be zero) for all $i$;
\item $|\la| = \sum_{i=1}^t |R_i|$.
\end{enumerate}
In this situation we introduce an alternative notation
\begin{equation*}
  \K_{\la;R}(q) = \K_{\la,\gamma,\eta}(q)
\end{equation*}
where $R$ stands for the sequence of partitions
$R=R(\eta,\gamma)=(R_1,R_2,\dots,R_t)$.
From now on we will use either notation as is convenient.
Say that $R$ is \textit{dominant} if $\gamma$ is.
\end{rem}

\subsection{Examples}
\label{special case section}

In each of the following examples, the sequence
of partitions $R$ consists entirely of \textit{rectangular} partitions.

\begin{enumerate}
\item (Kostka-Foulkes) Let $\gamma$ and $\la$ be partitions of $N$ of length
at most $n$ and $\eta=(1^n)$, so that the partition $R_i$
is a single row of length $\gamma_i$.  Then 
\begin{equation*}
  \K_{\la;R}(q) = K_{\la,\gamma}(q),
\end{equation*}
the Kostka-Foulkes polynomial (defined in \cite[III.6]{Mac}).
\item (Cocharge Kostka-Foulkes) Let $\eta$ be arbitrary, $\gamma=(1^n)$,
and $\la$ a partition of $n$,
so that $R_i = (1^{\eta_i})$ is a single column of length $\eta_i$.  Then 
\begin{equation*}
  \K_{\la;R}(q) = \tK_{\la^t,\eta^+}(q)
\end{equation*}
is the cocharge Kostka-Foulkes polynomial (defined in \cite[III.7]{Mac}),
where $\la^t$ is the transpose of the partition $\la$.
\item (Nilpotent orbit)  Let $\eta$ be arbitrary,
$k$ a positive integer,
$\gamma=(k^n)$, and $\la$ a partition of $kn$ with at most $n$ parts,
so that $R_i = (k^{\eta_i})$ is a rectangle with $\eta_i$ rows and $k$
columns.  Then $\K_{\la;R}(q)$ is the Poincar\'e polynomial of the
$(\la-k)$-th isotypic component of the coordinate ring
$\C[X_\mu]$ where $\mu=\eta^+$ \cite{We}.
\item (Two column Macdonald-Kostka) Let $2r\le n$.  Stembridge
\cite{St} (see also \cite{F}) showed that the two-column
Macdonald-Kostka polynomial can be written
\begin{equation*}
  K_{\la,(2^r,1^{n-2r})}(q,t) =
  	\sum_{k=0}^r q^k \begin{bmatrix} r\\ k \end{bmatrix}_t
	M^k_{r-k}(t)
\end{equation*}
where $\begin{bmatrix} r\\ k \end{bmatrix}_t$ is the usual $t$-binomial
coefficient and $M^d_m(t)$ is a polynomial in $t$ defined by
the recurrence
\begin{equation*}
\begin{split}
        M^0_m(t) &= K_{\lambda,(2^m,1^{n-2m})}(t) \\
        M^{d+1}_m(t) &= M^d_m(t) - t^{n-2m-d-1} M^d_{m+1}(t)
\end{split}
\end{equation*}
Using this defining recurrence for $M^d_m(t)$ it can be shown using
the methods of \cite{KKSW} that
\begin{equation}\label{small rects}
  M^d_m(t) = \K_{\la,((2)^m,(1,1)^d,(1)^{n-2m-2d})}(t),
\end{equation}
which involves only rectangles with at most two cells.  
The right hand side of \eqref{small rects}
satisfies the defining recurrence for $M^d_m(t)$, since
the following sequence of modules is exact, where
$\MM^d_m := \MM_{(1^m,2^d,1^{n-m-2d}),(1^m,0^{n-2m},(-1)^m)}$.
\begin{equation*}
\begin{split}
  0 \rightarrow 
  &\MM_{(1^m,2^d,1^{n-m-2d}),(1^m,0^{2d},1,0^{n-2d-2(m+1)},(-1)^{m+1})}
	[-(n-2m-1)] \rightarrow \\
  &\MM^d_m \rightarrow \MM^{d+1}_m \rightarrow 0.
\end{split}
\end{equation*}
The notation $[r]$ indicates a shift in homogeneous degree.
The proof is completed by establishing the graded character identity
\begin{equation*}
\begin{split}
  &H_{(1^m,2^d,1^{n-m-2d}),(1^m,0^{2d},1,0^{n-2d-2(m+1)},(-1)^{m+1})}(x;q)=\\
  &q^d H_{(1^{m+1},2^d,1^{n-(m+1)-2d}),(1^{m+1},0^{n-2(m+1)},(-1)^{m+1})}(x;q).
\end{split}
\end{equation*}
When $m=0$ the exact sequence is an explicit resolution of the ideal
of the nilpotent orbit closure $X_{(2^{d+1},1^{n-2(d+1)})}$
over the coordinate ring of the minimally larger one
$X_{(2^d,1^{n-2d})}$ \cite{KKSW}.
\end{enumerate}

\subsection{Positivity conjecture}

Broer has conjectured the following sufficient condition that
the polynomials $\K_{\la,\gamma,\eta}(q)$ have nonnegative
integer coefficients.

\begin{conj} \label{positivity conj} \cite{Broer2}
If $\gamma$ is dominant then in the notation of section \ref{Euler def},
\begin{equation*}
  H^i(G/B,\TTT' \otimes \LLL_\gamma)=0
\end{equation*}
for $i>0$.  In particular, the polynomial 
$\K_{\la,\gamma,\eta}(q)$, being the Poincar\'e polynomial
of an isotypic component of the graded module
$M_{\eta,\gamma}$, has nonnegative integer coefficients.
\end{conj}

This was verified by Broer in the case that
the vector bundle $\phi_* \LLL_\gamma$ is a line bundle \cite{Broer}.
In our case this means that each of the partitions $R_i$
is a rectangle.  We adopt a combinatorial approach to
positivity in section \ref{cat conj proof}.

\begin{ex}
\begin{enumerate}
\item Let $n=2$, $\la=(1,1)$, $\gamma=(0,2)$ and $\eta=(1,1)$.
Then $\K_{\la,\gamma,\eta}(q) = q-1$.
\item For $n=2$, $\la=(1,0)$, $\gamma=(0,1)$ and $\eta=(1,1)$,
$\K_{\la,\gamma,\eta}(q) = q$.
\item Let $n=3$, $\la=(2,1,0)$, $\gamma=(0,2,1)$ and $\eta=(1,1,1)$.
Then $\K_{\la,\gamma,\eta}(q) = q^3+q^2-q$.
\end{enumerate}
\end{ex}

\subsection{$q$-Kostant formula}

The following formula is a direct consequence
of formulas \eqref{gen func def}, \eqref{poly def},
and \eqref{Bott formula}.
Let $\epsilon_i$ be the $i$-th standard basis vector in $\Z^n$.

\begin{prop} \label{q Kostant}
\begin{equation*}
  \K_{\la,\gamma,\eta}(q) =
  \sum_{w\in W} (-1)^w
    \sum_{m:\Roots_\eta\rightarrow\N}
	q^{\sum_{(i,j)\in\Roots_\eta} m(i,j)}
\end{equation*}
where $m$ satisfies
\begin{equation*}
  \sum_{(i,j)\in\Roots_\eta} m(i,j)(\epsilon_i-\epsilon_j)=
	w^{-1}(\la+\rho)-(\gamma+\rho)
\end{equation*}
\end{prop}

In the Kostka-Foulkes special case, this is the formula for
Lusztig's $q$-analogue of weight multiplicity in type A
\cite[Ex. III.6.4]{Mac}.

\subsection{Generalized Morris recurrence}

We now derive a defining recurrence for the polynomials $\K_{\la;R}(q)$.
In the Kostka-Foulkes case this is due to Morris \cite{Mo}
and in the nilpotent orbit case, to Weyman \cite[(6.6)]{We}.
Following Remark \ref{normal parameters},
let us assume that each of the weights $R_i$ is a partition.

If $\eta=(\eta_1)$ consists of a single part,
then $\Roots_\eta$ is the empty set and $H_{\eta,\gamma}(x;q) =
s_\gamma(x)$.  In other words,
\begin{equation}
\label{base of recurrence}
\K_{\la;(R_1)}(q) = \delta_{\la,R_1},
\end{equation}
where $\delta$ is the Kronecker symbol.

Otherwise suppose that $\eta$ has more than one part.
Write $m = \eta_1$ and
\begin{equation*}
\begin{split}
  \etahat&=(\eta_2,\eta_3,\dots,\eta_t) \\
  \gammahat&=(\gamma_{m+1},\gamma_{m+2},\dots,\gamma_n) \\
  \Rhat&=(R_2,R_3,\dots,R_t)
\end{split}
\end{equation*}

For convenience let $y_i=x_i$ for $1\le i\le m$ and
$z_i=x_{m+i}$ for $1\le i\le n-m$.  Let $W_y$ and $W_z$
denote the subgroups of $W$ that act only on the $y$ variables
and $z$ variables respectively.  Let $\pi_x$, $\pi_y$, and $\pi_z$ be the
isobaric divided difference operators for the longest element
in $W$, $W_y$, and $W_z$ respectively.

Let $\la^*$ denote the highest weight of the contragredient
dual $(V_\la)^*$ of $V_\la$.  It is given explicitly by
$\la^* := w_0 (-\la) = (-\la_n,-\la_{n-1},\dots,-\la_1)$.  We have
\begin{equation*}
  s_{\la}(x^*) = s_{\la^*}(x)
\end{equation*}
where $x^*=(x_1^{-1},x_2^{-1},\dots,x_n^{-1})$.  Now
\begin{equation*}
\begin{split}
  x^\gamma B_\eta(x;q) &= y^{R_1} z^{\gammahat} \, B_{\etahat}(z;q) \,
    \prod_{1 \le i \le m < j \le n} (1-q\, x_i/x_j)^{-1} \\
    &= y^{R_1} z^{\gammahat} \, B_{\etahat}(z;q)
	\prod_{\substack{1 \le i \le m \\ 1\le j \le n-m}}
	(1-q \,y_i/z_j)^{-1} \\
  &= y^{R_1} z^{\gammahat} B_{\etahat}(z;q)
     \sum_\nu q^{|\nu|} s_\nu(y) s_{\nu^*}(z)
\end{split}
\end{equation*}
using the definition of $B_\eta(x;q)$ and Cauchy's formula.
The index variable $\nu$ runs over partitions of length at most
$\min(m,n-m)$.

Applying the operator $\pi_x = \pi_x \pi_y \pi_z$
to $x^\gamma B_\eta(x;q)$, we have 
\begin{equation}
\label{expansion}
\begin{split}
  H_{\eta,\gamma}(x;q) &= \pi_x \pi_y \pi_z 
	y^{R_1} z^{\gammahat} B_{\etahat}(z;q)
	\sum_\nu q^{|\nu|} s_\nu(y) s_{\nu^*}(z) \\
  &= \pi_x s_{R_1}(y) H_{\etahat,\gammahat}(z;q)
	\sum_\nu q^{|\nu|} s_\nu(y) s_{\nu^*}(z) \\
  &= \pi_x s_{R_1}(y) \sum_\nu q^{|\nu|} s_\nu(y) s_{\nu^*}(z)
     \sum_\sigma \K_{\sigma,\gammahat,\etahat}(q) s_\sigma(z) \\
  &= \pi_x \sum_\nu q^{|\nu|} \sum_\sigma
  \K_{\sigma,\gammahat,\etahat}(q)
  \sum_{\alpha,\beta} \LR^\alpha_{R_1,\nu} \LR^\beta_{\sigma\nu^*}
	s_\alpha(y) s_\beta(z)
\end{split}
\end{equation}
where $\alpha$ runs over the partitions of length at most $m$,
$\beta$ and $\sigma$ run over the dominant integral weights with
$n-m$ parts, and
\begin{equation} \label{LR}
	\LR^c_{ab} = \dim \Hom_{GL(n)}(V_c,V_a\otimes V_b).
\end{equation}
are the Littlewood-Richardson coefficients
for the dominant integral weights $a$, $b$, and $c$.

Taking the coefficient of $s_\la(x)$ on both sides of \eqref{expansion}
and applying \eqref{Bott formula} we have
\begin{equation}
\label{prerecurrence}
\begin{split}
  \K_{\la;R}(q) &=
  \sum_{w\in W/(W_y\times W_z)} (-1)^w q^{|\alpha(w)|-|R_1|}
  \sum_\sigma \K_{\sigma;\Rhat}(q) \\
  &\,\,\sum_{\substack{\nu \\ |\nu| = |\alpha(w)|-|R_1|}}
    \LR^{\alpha(w)}_{R_1,\nu} \LR^{\beta(w)}_{\sigma,\nu^*} 
\end{split}
\end{equation}
where $w$ runs over the minimal length coset representatives
and $\alpha(w)$ and $\beta(w)$ are the first $m$ and last $n-m$ parts
of the weight $w^{-1}(\la+\rho)-\rho$.
The restriction on $w$ is due to the fact that the formula
\eqref{Bott formula} is being applied only to Laurent
monomials of the form $y^\alpha z^\beta$,
where $\alpha$ and $\beta$ are dominant weights having
$m$ and $n-m$ parts respectively.  The $w$-th summand
is understood to be zero unless all the parts of $\alpha(w)$ are nonnegative
and $\alpha(w) \supseteq R_1$.  Note that $\beta(w)$ is always a partition
since $\la$ is.  Finally, the LR coefficients can be simplified.  Note that
\begin{equation*}
  \LR^c_{ab} = \dim \,(V_{c^*} \otimes V_a \otimes V_b)^{GL(n)}
\end{equation*}
where $(V_c)^*\cong V_{c^*}$ is the contragredient dual of $V_c$.
Applying duality and the definitions, one obtains
\begin{equation} \label{LR equalities}
  \LR^c_{ab} = \LR^c_{ba} = \LR^{c^*}_{a^*b^*} = \LR^a_{b^*c}.
\end{equation}
It follows that
\begin{equation*}
\begin{split}
  \sum_\nu \LR^{\alpha(w)}_{R_1,\nu} \LR^{\beta(w)}_{\sigma,\nu^*} 
&= \sum_\nu \LR^{\alpha(w)}_{R_1,\nu} \LR^{\sigma}_{\beta(w)\nu} \\
&= \dim \Hom_{GL(n-m)}(V_{\alpha(w)/R_1},V_{\sigma/\beta(w)}) \\
&= \dim \Hom_{GL(n-m)}(V_{\alpha(w)/R_1}\otimes V_{\beta(w)},V_\sigma) \\
&=: \LR^\sigma_{\alpha(w)/R_1,\beta(w)}
\end{split}
\end{equation*}
where $V_{\la/\mu}$ is the $GL(n-m)$-module whose formal character
is the skew Schur polynomial $s_{\la/\mu}(z)$.
Then \eqref{prerecurrence} can be expressed as
\begin{equation}
\label{recurrence}
  \K_{\la;R}(q) =
  \sum_{w\in W/(W_y\times W_z)} (-1)^w q^{|\alpha(w)|-|R_1|}
  \sum_\sigma \LR^\sigma_{\alpha(w)/R_1,\beta(w)} \K_{\sigma;\Rhat}(q)
\end{equation}
Clearly the recurrence \eqref{recurrence} together with the initial
condition \eqref{base of recurrence} uniquely defines the polynomials
$\K_{\la;R}(q)$.

\subsection{$q$-analogue of LR coefficients}

Let $\inner{\,}{\,}$ denote the Hall inner product on
the ring $\C[x]^W$ of $W$-symmetric polynomials.
We employ the notation of Remark \ref{normal parameters}.
The following result does not assume that $R$ is dominant.

\begin{prop} \label{LR specialization} 
\begin{equation} \label{shape product}
  \K_{\la;R}(1) = \LR^\la_R :=
	\inner{s_\la(x)}{s_{R_1}(x) s_{R_2}(x) \dots s_{R_t}(x)}
\end{equation}
\end{prop}
\begin{proof} The proof proceeds by induction on $t$, the number of
partitions in $R$.  For $t=1$ the result holds by \eqref{base of recurrence}.
Suppose that $t>1$.  By \eqref{recurrence} and induction we have
\begin{equation*}
\begin{split}
\K_{\la;R}(1) &=
  \sum_{w\in W/(W_y\times W_z)}
    (-1)^w \sum_\sigma \K_{\sigma;\Rhat}(1) \,
	\LR^\sigma_{\alpha(w)/R_1,\beta(w)} \\
 &= \sum_{w\in W/(W_y\times W_z)}
   (-1)^w \sum_\sigma \inner{s_\sigma}{s_{R_2} s_{R_3} \dots s_{R_t}}
	\inner{s_\sigma}{s_{\alpha(w)/R_1} s_\beta(w)} \\
 &= \inner{s_{R_2} s_{R_3} \dots s_{R_t}}{
	\sum_{w\in W/(W_y\times W_z)} (-1)^w s_{\alpha(w)/R_1}\, s_\beta(w)}
\end{split}
\end{equation*}
The Jacobi-Trudi formula for the skew Schur polynomial
$s_{\la/\mu}$ is given by the determinant
\begin{equation*}
  s_{\la/\mu} = \det h_{\la_i-i-(\mu_j-j)}(x)
\end{equation*}
Using Laplace's expansion in terms of $m \times m$ minors
involving the first $m$ columns, we have
\begin{equation*}
  s_{\la/R_1} =
  \sum_{w\in W/(W_y \times W_z)}
  (-1)^w s_{\alpha(w)/R_1} s_\beta(w)
\end{equation*}
from which \eqref{shape product} follows.
\end{proof}

\begin{rem} Proposition \ref{LR specialization} can be
proven using the definition of the module $\MM_{\eta,\gamma}$
together with an algebraic reciprocity theorem.
\end{rem}

\subsection{Contragredient duality}

The following symmetry of the Poincar\'e polynomial
$\K_{\la,\gamma,\eta}$ is an immediate consequence of the formulas
\eqref{gen func def} and \eqref{poly def}.

\begin{prop} \label{weight duality}
\begin{equation}
  \K_{\la,\gamma,\eta}(q) =
  \K_{\la^*,\gamma^*,(\eta_t,\dots,\eta_1)}(q).
\end{equation}
\end{prop}

This gives a $q$-analogue of the ``box complement" duality of LR
coefficients.  Suppose that $\la$ and $\gamma$ are partitions and
$R=R(\eta,\gamma)$ is the associated sequence of partitions.
Let $M$ be a positive integer such that $M\ge \max(\la_1,\gamma_1)$.  
Let $\widetilde{\la}$ be the partition obtained by taking the complement
of the diagram of $\la$ inside the $n\times M$ rectangle
and rotating it 180 degrees.  That is, $\widetilde{\la}=\la^*+M$ in the
notation of section \ref{normal parameters}.  Define 
$\widetilde{R_i}$ similarly, except that the latter is complemented
inside the $\eta_i\times M$ rectangle.  Then
\begin{equation} \label{box complement}
  \K_{\la;R}(q) = \K_{\widetilde{\la};
	(\widetilde{R_t},\dots,\widetilde{R_1})}(q)
\end{equation}

\subsection{Symmetry}

\begin{prop} \label{symm prop}
Let $R=R(\eta,\gamma)$ be dominant
and $R'$ a dominant reordering of $R$.  Then for any $\la$,
\begin{equation*}
  \K_{\la;R}(q) = \K_{\la;R'}(q).
\end{equation*}
\end{prop}
\begin{proof} From the assumptions it follows that $R'$ can be written
$R'=R(\eta',\gamma)$.  The dominance condition implies
that $R'$ is obtained from $R$ by a sequence of exchanges
of adjacent rectangular partitions that have the same number of columns.
A series of reductions can be made.  
It may be assumed that $R$ and $R'$ differ by one such exchange,
so that $R'$ is obtained from $R$ by exchanging $R_i$ and $R_{i+1}$, say.
If $i>1$ then by the generalized Morris-Weyman recurrence
\eqref{recurrence}, the first common partition $R_1=R_1'$ may be removed,
so that by induction it may be assumed that $i=1$.
On the other hand, by applying the ``box complement"
formula \eqref{box complement} it may be assumed that $R$ and $R'$
differ by exchanging their last two partitions.  As before the common
partitions at the beginning may be removed, so that one may assume that
$R=(R_1,R_2)$ and $R'=(R_2,R_1)$ where $R_1$ and $R_2$ are rectangles
having the same number of columns.

Without loss of generality it may be assumed that $\eta_1>\eta_2$,
so that $\eta=(\eta_1,\eta_2)$ is a partition and $\eta'=(\eta_2,\eta_1)$.
By \eqref{shift all} it is enough to show that
\begin{equation*}
  H_{\eta,(0,0)}(x;q) = H_{\eta',(0,0)}(x;q),
\end{equation*}
These are the formal characters of the Euler characteristics
of the structure sheafs of the desingularizations $Z(\eta)$ and
$Z(\eta')$ of $X_\eta$.  By \cite{Broer} the higher direct images of
the corresponding desingularization maps (call them $q$ and $q'$)
vanish, so
\begin{equation*}
  \chi_{\eta,(0,0)} = q_*(\OOO_{Z_\eta}) = \OOO_{X_\eta} =
  	q'_*(\OOO_{Z_{\eta'}}) = \chi_{\eta',(0,0)},
\end{equation*}
and we are done.

The case that $R$ has two partitions, may also be verified by the explicit
formula \eqref{two part Poincare} derived below.
\end{proof}

This result is not at all obvious when looking directly at the
definition of $H_{\eta,\gamma}(x;q)$ and $H_{\eta',\gamma}(x;q)$,
since the sets of weights $\Roots_\eta$ and $\Roots_{\eta'}$ differ.

\subsection{Monotonicity}

If $\mu$ and $\nu$ are partitions such that $\mu\dom\nu$
(that is, $\mu_1+\dots+\mu_i\ge\nu_1+\dots+\nu_i$ for all $i$)
then $X_\mu \subset X_\nu$, so that restriction of functions gives
a natural graded $G$-module epimorphism $\C[X_\nu]\rightarrow\C[X_\mu]$.
There are similar epimorphisms for the twisted modules
$M_{\eta,\gamma}$.  These yield inequalities for the Poincar\'e polynomials
of their isotypic components.

Unfortunately the Poincar\'e polynomial $\K_{\la,\gamma,\eta}(q)$
comes from the Euler characteristic $\chi_{\eta,\gamma}$ and not
not just the term $M_{\eta,\gamma}$ in cohomological degree zero.
So whenever Conjecture \ref{positivity conj} holds
(for example, for dominant sequences of rectangles),
then one has a corresponding inequality for the polynomials
$\K_{\la,\gamma,\eta}(q)$.  Here are two such inequalities.

\begin{conj} Let $R=R(\eta,\gamma)$ be dominant and
$R'=R(\eta',\gamma)$, where
$\eta'$ is obtained from $\eta$ by replacing some part $\eta_i$
by another sequence that sums to $\eta_i$.  Then for any $\la$,
\begin{equation*}
 \K_{\la;R}(q) \le \K_{\la;R'}(q)
\end{equation*}
coefficientwise.
\end{conj}

\begin{conj} \label{second mono}
Suppose that $R=R(\eta,\gamma)$ is dominant
and $R'$ is obtained from $R$ by replacing some subsequence
of the form 
\begin{equation*}
	((k^{\alpha_1}),\dots,(k^{\alpha_l}))
\end{equation*}
by
\begin{equation*}
	((k^{\beta_1}),\dots,(k^{\beta_l}))
\end{equation*}
where $|\alpha|=|\beta|$ and $\alpha^+ \dom \beta^+$.  Then for any $\la$,
\begin{equation*}
 \K_{\la;R}(q) \le \K_{\la;R'}(q)
\end{equation*}
coefficientwise.
\end{conj}

In the cocharge Kostka-Foulkes case this is a well-known monotonicity
property for Kostka-Foulkes polynomials \eqref{monotonicity} that will be
discussed at length in section \ref{app sec}.

\section{Combinatorics}

In this section we give a conjectural combinatorial description
(Conjecture \ref{cat conj}) of the polynomials $\K_{\la;R}(q)$
as well as a general approach for its proof, which is shown to
succeed in special cases.  These require considerable combinatorial
preliminaries.  The new material begins in subsection \ref{cat sec}.

\subsection{Tableaux and RSK}

We adopt the English convention for partition diagrams and tableaux.
The \textit{Ferrers diagram} of a partition
$\la=(\la_1\ge\la_2\ge\dots\ge\la_n\ge0)$
is the set of ordered pairs of integers
$D(\la)=\{(i,j): 1\le j\le \la_i\}$.  A \textit{skew shape}
$\la/\mu$ is the set difference $D(\la)-D(\mu)$
of Ferrers diagrams of partitions.
A \textit{(skew) tableau} $T$ is a function $T:D\rightarrow\N_+$ from
a skew shape $D$ to the positive integers.  The domain of a tableau $T$
is called its \textit{shape}.
A tableau $T$ of shape $D$ is depicted as a partial matrix whose
$(i,j)$-th position contains the value $T(i,j)$ for $(i,j)\in D$.
A tableau is \textit{column strict} if it weakly increases from left to
right within each row and strictly increases from top to bottom within each
column.  The \textit{row-reading word} of a (skew) tableau $T$ is
$\dots w^2 w^1$, where $w^i$ is the word comprising the $i$-th row of $T$,
read from left to right.  The \textit{content} of a word $w$ (or tableau
$T$) is the sequence $(c_1,c_2,\dots)$ where $c_i$ is the number of
occurrences of the letter $i$ in $w$ (or $T$).  Say that a word of length
$n$ is \textit{standard} if it has content $(1^n)$.
A tableau whose shape consists of $n$ cells
is said to be \textit{standard} if it is column strict and
has content $(1^n)$.

\begin{ex} Let $D=\la/\mu$ where $\la=(6,5,3,3)$ and $\mu=(3,2)$.
A column strict tableau $T$ of shape $D$ is depicted below.
\begin{equation*}
T = \begin{matrix}
\sq & \sq & \sq & 1 & 2 & 2 \\
\sq & \sq & 1 & 2 & 3 &   \\
2 & 3 & 3 &   &   &   \\
4 & 4 & 5 &   &   & 
\end{matrix}
\end{equation*}
The row-reading word of $T$ is
\begin{equation*}
  445\,233\,123\,122
\end{equation*}
and the content of $T$ is $(2,4,3,2,1)$.
\end{ex}

The \textit{Knuth equivalence} is the equivalence relation $\Kn$ on words
that is generated by the relations of the following form,
where $u$ and $v$ are arbitrary words and $x$, $y$, and $z$ are letters:
\begin{equation*}
\begin{split}
  u x z y v &\Kn u z x y v \qquad\text{for $x\le y<z$} \\
  u y x z v &\Kn u y z x v \qquad\text{for $x<y\le z$}
\end{split}
\end{equation*}

For the word $u$, let $P(u)$ be Schensted's $P$-tableau, that is,
the unique column strict tableau of partition shape whose row-reading
word is Knuth equivalent to the word $u$.

We establish some notation for the column insertion version of the
Robinson-Schensted-Knuth (RSK) correspondence.
For each $i\ge 1$, let $u^i$ be a weakly increasing word (almost all empty).
The column insertion RSK correspondence is the bijection
from the set of such sequences of words,
to pairs of column strict tableaux $(P,Q)$ of the same shape,
defined by $P = P(\dots u^2 u^1)$ and
$\shape(Q|_{[i]}) = \shape(P(u^i u^{i-1} \dots u^1))$ for all $i$,
where $Q|_{[i]}$ denotes the restriction of the tableau $Q$ to
the letters in the set $[i]$.
In other words, to produce $P$ one performs the column insertion of
the word $\dots u^2 u^1$ starting from the right end, and to
produce $Q$ one records the insertions of letters in the subword $u^i$
by the letter $i$ in $Q$.  In particular, if $\alpha_i$ is the length
$|u^i|$ of the word $u^i$, then $Q$ has content $\alpha$.

\subsection{Crystal operators}

We recall the definitions of the $r$-th crystal reflection, raising,
and lowering operators $s_r$, $e_r$, and $f_r$ on words.
These are due to Lascoux and Sch\"utzenberger \cite{LS2} \cite{LLT2}.

Let $r$ be a positive integer and $u$ a word.  Ignore all letters
of $u$ which are not in the set $\{r,\rp \}$.  View each occurrence
of the letter $r$ (resp. $r+1$) as a right (resp. left) parenthesis.
Perform the usual matching of parentheses.
Say that an occurrence of a letter $r$ or $\rp$
in $u$ is \textit{$r$-paired} if it corresponds to a matched
parenthesis.  Otherwise call
that letter \textit{$r$-unpaired}.  It is easy to see that the subword
of $r$-unpaired letters of $u$ has the form $r^p \rp^q$ where
$r^p$ denotes the word consisting of $p$ occurrences of the letter $r$.

Consider the three operators $s_r$, $e_r$, and $f_r$ on words,
which are called the $r$-th crystal reflection,
raising, and lowering operators respectively.
Each is applied to the word $u$ by replacing the 
$r$-unpaired subword $r^p \rp^q$ of $u$ by another subword of the
same form.  Below each operator is listed, together with the subword
that it uses to replace $r^p \rp^q$ in $u$.
\begin{enumerate}
\item For $s_r u$ use $r^q \rp^p$.  This operator clearly switches
the number of $r$'s and $\rp$'s in $u$.
\item For $e_r u$, use $r^{p+1} \rp^{q-1}$.  This
only makes sense if $q>0$, that is, there is an $r$-unpaired letter
$\rp$ in $u$.
\item For $f_r u$, use $r^{p-1} \rp^{q+1}$.  This makes sense
when $p>0$, that is, there is an $r$-unpaired $r$ in $u$.
\end{enumerate}

\begin{ex} For $r=2$, the $r$-th crystal operators are calculated on the
word $u$.  The $r$-unpaired letters are underlined.
\begin{equation}
\begin{split}
u&=1\ul 24312\ul 2\ul 3342\ul 3\ul 34\ul 3313123422\ul 3 \\
s_2\, u&=1\ul 24312\ul 2\ul 2342\ul 2\ul 24\ul 3313123422\ul 3 \\
e_2\, u&=1\ul 24312\ul 2\ul 2342\ul 3\ul 34\ul 3313123422\ul 3 \\
f_2\, u&=1\ul 24312\ul 3\ul 3342\ul 3\ul 34\ul 3313123422\ul 3
\end{split}
\end{equation}
\end{ex}

These operators may be defined on skew column strict tableaux by
acting on the row-reading word.  Each produces a column strict tableau
of the same shape as the original (skew) tableau.  It is proven in
\cite{LS2} that $\{s_1,s_2,\dots\}$ satisfy the
Moore-Coxeter relations as operators on words, so that one may
define an action of the infinite symmetric group on words by
\begin{equation} \label{plactic action}
  w u = s_{i_1} s_{i_2} \dots s_{i_p} u
\end{equation}
where $u$ is a word, $w$ is a permutation and
$w=s_{i_1}\dots s_{i_p}$ is a reduced decomposition of $w$.

We say that two words (or tableaux) are in the same
\textit{$r$-string} if
each can be obtained from the other by a power of $e_r$ or $f_r$,
or equivalently, they differ only in their $r$-unpaired subwords.

\subsection{Lattice property}

Let $\mu$ be a partition.  Say that the word $u$ is  \textit{$\mu$-lattice}
if the sum of $\mu$ and the content of every final subword of $u$
is a partition.  Say that a (skew) column strict tableau
is $\mu$-lattice if its row-reading word is.  Say that a word or tableau
is \textit{lattice} if it is $\mu$-lattice when $\mu$ is the empty partition.

The following formulation of the Littlewood-Richardson rule is well-known.
See \cite{FG}, \cite{Sh} and \cite{Z}.

\begin{thm} (LR rule) The LR coefficient
\begin{equation*}
  \inner{\la/\mu}{\sigma/\tau}
\end{equation*}
is given by the number of column strict tableaux of shape $\sigma/\tau$
of content $\la-\mu$ that are $\mu$-lattice.
\end{thm}

The following is a reformulation of a result of D. White \cite{White}.

\begin{thm} \label{fitting}
Let $\{u^i\}$ be a sequence of weakly increasing words and
$(P,Q)$ the corresponding tableau pair under column RSK.
Then there is a column-strict tableau of shape $\la/\mu$,
whose $i$-th row is given by the word $u^i$ for all $i$,
if and only if the content of $Q$ is $\la-\mu$ and
$Q$ is $\mu$-lattice.
\end{thm}

The proof of the following result is straightforward and left to
the reader; see also \cite{RS} and \cite{SW}.
 
\begin{lem} \label{lattice}
\begin{enumerate}
\item A word is $\mu$-lattice if and only if the number of $r$-unpaired
$\rp$'s is at most $\mu_r - \mu_{\rp}$ for all $r$.
\item There is an involution on the set of words that are not $\mu$-lattice,
given by $u\mapsto s_r e_r^{\mu_r-\mu_{\rp}+1} u$, where $r+1$ is the
rightmost letter in $u$ where $\mu$-latticeness fails.
\end{enumerate}
\end{lem}

\subsection{Evacuation}

Let $T$ be a column strict tableau of 
content $\alpha=(\alpha_1,\alpha_2,\dots,\alpha_n)$
and partition shape.
The evacuation $\ev_{[n]}(T)$ of $T$ with respect to the alphabet
$[n]$ is the unique column strict tableau in the alphabet $[n]$
such that
\begin{equation*}
  \shape((\ev_{[n]}(T))|_{[i]}) = \shape(P(T|_{[n+1-i,n]})).
\end{equation*}
Clearly $\ev(T)$ has content $(\alpha_n,\dots,\alpha_1)$
and the same shape as $T$.

The main result on evacuation is the following.

\begin{thm} \label{ev thm} \cite{LS2}
Let $\{u^i: 1\le i \le n \}$ be a collection of weakly increasing words
in the alphabet $[N]$ and $(P,Q)$ the corresponding pair of tableaux under
column RSK.  Let $v^i$ be the reverse of the word obtained from
$u^{n+1-i}$ by complementing each letter in the alphabet $[N]$.
Then the sequence of words $\{v^i\}$ corresponds under column RSK
to the tableau pair $(\ev_{[N]}(P),\ev_{[n]}(Q))$.
\end{thm}

\subsection{Two row jeux-de-taquin}
\label{two row stuff}

There is a duality between the crystal operators and
jeux-de-taquin on two-row skew column strict tableaux.
This is described below.

Define the \textit{overlap} of the pair $(v,u)$ of weakly increasing
words to be the length of the second row in the tableau $P(v u)$,
or equivalently, the maximum number of columns of size two among the
skew column strict tableaux with first row $u$ and second row $v$.

\begin{lem} \label{overlap} Let $(P,Q)$ be the tableau pair
obtained by column RSK from the sequence of words
$\{v^i\}$.  Then the overlap of the pair of words
$(v^{r+1},v^r)$ is equal to the number of $r$-pairs in $Q$.
\end{lem}
\begin{proof} (Sketch) In the case $r=1$ this is easy to check
directly.  So suppose $r>1$.  Define a new sequence of weakly
increasing words $\{u^i\}$ by $u^i=v^i$ for $i\ge r$ and
letting $u^i$ be empty for $i < r$.  Let $(\widehat{P},\widehat{Q})$ be
the resulting tableau pair.  Then it can be shown \cite{H} \cite{Sh}
that $\widehat{Q} = P(Q|_{[r,n]})$, that is, $\widehat{Q}$ is obtained from
$Q$ by removing all letters strictly less than $r$ and then taking the
Schensted $P$-tableau.  Moreover it is straightforward
to show that the number of $r$-pairs is invariant under Knuth equivalence;
this is equivalent to the well-known fact (see \cite{Sa}) that the lattice
property is invariant under Knuth equivalence.
This reduces the proof to the case $r=1$.
\end{proof}

\begin{lem} \label{two row dual} Let $\{v^i\}$ and $(P,Q)$ be
as in Lemma \ref{overlap} and let $\{{v'}^i\}$ be another
sequence of weakly increasing words with corresponding tableau
pair $(P',Q')$.  The following are equivalent.
\begin{enumerate}
\item $P=P'$, and $Q$ and $Q'$ are in the same $r$-string.
\item ${v'}^i = v^i$ for $i\not\in[r,r+1]$ and
$P({v'}^{r+1} {v'}^r) = P(v^{r+1} v^r)$.
\end{enumerate}
\end{lem}
\begin{proof} (Sketch) Immediately one may reduce to the
case that $r+1$ is the largest letter in $Q$.
If $r=1$ then the proof is trivial,
since any two column strict tableaux
of the same partition shape in the alphabet $[1,2]$ are in the
same $1$-string.  Otherwise suppose $r>1$.  Let $N$ be the largest
letter of $P$ and $u^i$ the weakly increasing word given by the
reverse of the complement (in the interval $[N]$) of the word $v^{r+2-i}$
for $1\le i\le r+1$.  By Theorem \ref{ev thm}
the corresponding tableau pair under column RSK is
given by $(\ev_{[N]}(P),\ev_{[r+1]}(Q))$.  It can be shown (see the proof
of Lemma 63 in \cite{RnS}) that 
\begin{equation*}
\begin{split}
  \ev_{[r+1]} \circ e_i &= f_{r+1-i} \circ \ev_{[r+1]} \\
  \ev_{[r+1]} \circ f_i &= e_{r+1-i} \circ \ev_{[r+1]}
\end{split}
\end{equation*}
for any $1\le i\le r$.  Setting $i=1$,
the proof may be reduced to the case $r=1$.
\end{proof}

\subsection{Charge}

Following \cite{LS2} we review the definition of the
\textit{charge}, an $\N$-valued function on words.
There are three parts to the definition: words of standard content,
partition content, and arbitrary content.

Suppose first that $u$ is a word of content $(1^n)$.
Affix an index $c_i$ to the letter $i$ in $u$ according to the
rule that $c_1=0$ and $c_i=c_{i-1}$ if $i$ appears to the left of $i-1$ in
$u$ and $c_i=c_{i-1}+1$ if $i$ appears to the right of $i-1$ in $u$.  Let
\begin{equation*}
	\charge(u) = c_1 + c_2 + \dots + c_n
\end{equation*}

If $u$ has partition content $\mu$, then define
\begin{equation*}
	\charge(u) = \charge(u^1)+\charge(u^2)+\dotsm
\end{equation*}
where $u$ is partitioned into disjoint standard subwords
$u^j$ of length $\mu^t_j$ using the following left circular reading.
To compute $u^1$, start from the right end of $u$ and scan to the
left.  Choose the first $1$ encountered, then the first
$2$ that occurs to the left of the selected letter $1$, etc. 
If at any point there is no $i+1$ to the left of the selected
letter $i$, circle around to the right end of $u$ and continue
scanning to the left.  This process selects the subword $u^1$
of $u$.  Erase the letters of $u^1$ from $u$ and repeat this
process, obtaining the subword $u^2$.  Continue
until all the letters of $u$ have been exhausted.  

\setcounter{MaxMatrixCols}{25}

\begin{ex} The charge is calculated on the word $u$.
The words $u$, $u^1$, and $u^2$ appear below.
\begin{equation*}
\begin{matrix}
u  &=&4&3&2&3&4&1&1&2&5&5 \\
u^1&=&4&3&2& & & &1& & &5 \\
u^2&=& & & &3&4&1& &2&5&
\end{matrix}
\end{equation*}

\setcounter{MaxMatrixCols}{15}

The charges of the subwords $u^1$ and $u^2$ are calculated.
Each index $c_i$ is written below the letter $i$.
\begin{equation*}
\begin{matrix}
4&3&2&1&5\\
0&0&0&0&1
\end{matrix}
\qquad
\begin{matrix}
3&4&1&2&5\\
1&2&0&1&3
\end{matrix}
\end{equation*}
So
\begin{equation*}
  \charge(u)=\charge(u^1)+\charge(u^2) = 1 + 7 = 8.
\end{equation*}
\end{ex}

Finally, if the word $u$ has content $\alpha$,
define $\charge(u) = \charge((w_\alpha)^{-1} u)$ where the permutation
$(w_\alpha)^{-1}$ acts on the word $u$ as in \eqref{plactic action}.

\begin{thm} \label{charge thm} \cite{LS1}
\begin{equation*}
  K_{\la,\mu}(q) = \sum_T q^{\charge(T)}
\end{equation*}
where $T$ runs over the set of column strict tableaux
of shape $\la$ and content $\mu$.
\end{thm}

The charge has the following intrinsic characterization.

\begin{thm} \label{charge characterization} \cite{LS2}
The charge is the unique function from 
words to $\N$ such that:
\begin{enumerate}
\item For any word $u$ and any permutation $w$,
$\charge(w u) = \charge(u)$ where $w u$ is defined in
\eqref{plactic action}.
\item The charge of the empty word is zero.
\item If $u$ is a word of partition content $\mu$ of the form
$v 1^{\mu_1}$, then 
\begin{equation*}
  \charge(u) = \charge(v)
\end{equation*}
where $v$ is regarded
as a word of partition content $(\mu_2,\mu_3,\dots)$ in the alphabet
$\{2,3,\dots\}$.
\item Let $a>1$ be a letter and $x$ a word such that the word $ax$ has
partition content.  Then 
\begin{equation*}
  \charge(xa)=\charge(ax)+1.
\end{equation*}
\item The charge is constant on Knuth equivalence classes.
\end{enumerate}
\end{thm}

Define the charge of a (skew) column strict tableau to be
the charge of its row-reading word.

\subsection{Catabolizable tableaux}
\label{cat sec}

We now generalize the definition of a catabolizable tableau
given in \cite{SW}, which was inspired by the catabolism construction
of Lascoux and Sch\"utzenberger \cite{La} \cite{LS2}.

Let $R=(R_1,R_2,\dots,R_t)$ be a sequence of partitions determined
by the pair $(\eta,\gamma)$ as in Remark \ref{normal parameters}.
Let $Y_i$ be the tableau of shape $R_i$
whose $j$-th row is filled with the $j$-th largest letter of the
subinterval $A_i$.  As before we write $\eta=(m,\etahat)$.

\begin{ex} \label{R seq ex}
Let $n=5$, $\eta=(2,2,1)$, $\gamma=(3,2,2,1,1)$, so that
$m=2$ and $R=((3,2),(2,1),(1))$.  Then $A_1=[1,2]$, $A_2=[3,4]$, and
$A_3=[5,5]$.  The tableaux $Y_i$ are given by
\begin{equation*}
  Y_1 = \begin{matrix} 1&1&1\\2&2& \end{matrix} \qquad
  Y_2 = \begin{matrix} 3&3\\4& \end{matrix} \qquad
  Y_3 = \begin{matrix} 5 \\ \end{matrix}
\end{equation*}
\end{ex}

Given a (possibly skew) column strict tableau $T$ and index $r$,
let $H_r(T) = P(T_n T_s)$ where $T_n$ and $T_s$ are the
north and south subtableaux obtained by slicing $T$
horizontally between its $r$-th and $(r+1)$-st rows.

Let $S$ be a column strict tableau of partition shape in the alphabet $[n]$.
Suppose the restriction $S|_{A_1}$ of $S$ to the subalphabet $A_1$
is the tableau $Y_1$.  In thise case the $R_1$-catabolism of $S$
is defined to be the tableau $\cat_{R_1}(S)=H_m(S-Y_1)$.

The notion of a $R$-catabolizable tableau is uniquely defined
by the following rules.
\begin{enumerate}
\item If $R$ is the empty sequence, then the unique
$R$-catabolizable tableau is the empty tableau.
\item Otherwise, $T$ is $R$-catabolizable if and only if
$T|_{[m]} = Y_1$ and $\cat_{R_1}(T)$ is
$\Rhat$-catabolizable in the alphabet $[m+1,n]$.
\end{enumerate}

Denote by $CT(\la;R)$ the set of $R$-catabolizable tableaux of
shape $\la$.

\begin{ex} Continuing the previous example, let $\la=(5,3,1,0,0)$.
The four $R$-catabolizable tableaux of shape $\la$ are:
\begin{equation*}
\begin{matrix}
  1&1&1&3&4\\
  2&2&5& & \\
  3& & & &
\end{matrix}
\qquad
\begin{matrix}
  1&1&1&3&3\\
  2&2&4& & \\
  5& & & &
\end{matrix}
\qquad
\begin{matrix}
  1&1&1&4&5\\
  2&2&3& & \\
  3& & & &
\end{matrix}
\qquad
\begin{matrix}
  1&1&1&3&5\\
  2&2&4& & \\
  3& & & &
\end{matrix}
\end{equation*}
Let $S$ be the last tableau.  It is shown to be $R$-catabolizable
as follows.
\begin{equation*}
  S = \begin{matrix}
  1&1&1&3&5\\
  2&2&4& & \\
  3& & & &
\end{matrix}
\qquad
  S_n = \begin{matrix}
  \sq&\sq&\sq&3&5\\
  \sq&\sq&4& & \\
   & & & &
\end{matrix}
\qquad
  S_s = \begin{matrix}
  \sq&\sq&\sq& & \\
  \sq&\sq& & & \\
  3& & & &
\end{matrix}
\end{equation*}
\begin{equation*}
\cat_{R_1}(S) = P(435\,\,3) =
\begin{matrix}
3&3\\
4&5
\end{matrix} 
\end{equation*}
Now $\cat_{R_1}(S)$ contains $Y_2$ and
and $\cat_{R_2} \cat_{R_1}(S) = 5$.  The latter tableau 
contains $Y_3$ and
$\cat_{R_3} \cat_{R_2} \cat_{R_1}(S)$ is the empty tableau.
Thus $S$ is $R$-catabolizable.
\end{ex}

It is clear that any $R$-catabolizable tableau has content $\gamma$.

\begin{conj} \label{cat conj} For $R$ dominant,
\begin{equation*}
	\K_{\la;R}(q) = \sum_S q^{\charge(S)}
\end{equation*}
where $S$ runs over the set of $R$-catabolizable tableaux
of shape $\la$.
\end{conj}

\begin{ex} For the previous example,
$\K_{\la;R}(q) = q^3 + 3 q^4$.  This can be verified by computing
the polynomial two ways: 1) using the recurrence \eqref{recurrence}
and 2) evaluating the charge statistic on the four tableaux above.
The charges of the tableaux are (in order) $3,4,4,4$.
\end{ex}

We prove Conjecture \ref{cat conj} in the following cases.
\begin{enumerate}
\item $\eta$ is a hook partition, that is, $\eta_i=1$ for $i>1$.
An important subcase is the Kostka-Foulkes case, where $\eta=(1^n)$ and
$R_i$ is a single row of length $\gamma_i$.  Then an
$R$-catabolizable tableau is simply a column strict tableau of content
$\gamma$, and Conjecture \ref{cat conj} reduces to a
theorem of Lascoux and Sch\"utzenberger \cite{LS1}.
Even in this special case our proof differs from that of Lascoux and
Sch\"utzenberger, whose gaps were bridged by Butler \cite{Bu}.
\item In the cocharge Kostka-Foulkes case, where $R_i$
is a single column of length $\eta_i$ and $\eta$ is a partition,
Conjecture \ref{cat conj} reduces to a formula that is very nearly the
same as one of Lascoux \cite{La} for the cocharge Kostka-Foulkes
polynomials.  However the equivalence of the combinatorial definitions of
the corresponding sets of standard tableaux can only be resolved with
considerable effort.  This is done in section \ref{cyc std proof}.
\item $\eta$ has two parts.
\end{enumerate}

\subsection{Proof strategy for Conjecture \ref{cat conj}}
\label{cat conj proof}

We give a general approach for a proof of Conjecture
\ref{cat conj} using a sign-reversing involution.  For this
purpose it is convenient to modify the recurrence \eqref{recurrence}
to include more terms.  The complete expansion of the Jacobi-Trudi
determinant for the skew Schur polynomial $s_{\la/R_1}$ yields
\begin{equation}
\label{new rec}
  \K_{\la;R}(q) = \sum_{w\in W} (-1)^w q^{|\alpha(w)|-|R_1|}
  	\sum_\sigma K_{\sigma,(\alpha(w)-R_1,\beta(w))}
  		\K_{\sigma;\Rhat}(q)
\end{equation}
where $\alpha(w)$ and $\beta(w)$ are the first $m$ and last
$n-m$ parts of the weight
\begin{equation*}
  \xi(w) := w^{-1}(\la+\rho)-\rho
\end{equation*}
and $K_{\la,\alpha}$ is the Kostka number \cite[I.6]{Mac},
the number of column strict tableaux of shape $\la$ and content $\alpha$.

We now combinatorialize both sides of \eqref{new rec}.
Let $\SSS'$ be the set of triples $(w,T,U)$ where $w\in W$
and $T$ and $U$ are column strict tableaux of the same partition shape
where $T$ has content $(0^m,\gammahat)$ and $U$ has content
$(\beta(w),\alpha(w)-R_1)$.  Let $\SSS \subset \SSS'$ be the
subset of triples $(w,T,U)$ such that $T$ is $\Rhat$-catabolizable
in the alphabet $[m+1,n]$.  The reordering of parts of the content of $U$
is justified since the Kostka number is symmetric in its second index.
Define a sign and weight on $\SSS'$ by
\begin{equation*}
\begin{split}
  \sign(w,T,U) &= (-1)^w \\
  \weight(w,T,U) &= q^{|\alpha(w)|-|R_1| + \charge(T)}
\end{split}
\end{equation*}
where $T$ is regarded as a tableau of partition content in the
alphabet $[m+1,n]$.
By induction the right hand side of \eqref{new rec} is given by
\begin{equation*}
  \sum_{(w,T,U)\in\SSS} \sign(w,T,U) \weight(w,T,U)
\end{equation*}

\begin{ex} Let $n=8$, $\eta=(2,2,2,1,1)$, and $\gamma=(3^8)$, so that
$m=2$ and
\begin{equation*}
\begin{split}
  R&=((3,3),(3,3),(3,3),(3),(3)) \\
  \Rhat&=((3,3),(3,3),(3),(3)).
\end{split}
\end{equation*}

Let $w=32154678$ and $\la=(6,5,5,5,2,1,0,0)$\, so that
$\xi(w) = (3,5,8,1,6,1,0,0)$, $\alpha(w)=(3,5)$ and
$\beta(w)=(8,1,6,1,0,0)$.  Let $T$ and $U$ be given by
\begin{equation*}
  T = \begin{matrix}
    3&3&3&5&6&7&7&7\\
    4&4&4&6&8&8& & \\
    5&5&8& & & & & \\
    6& & & & & & & 
  \end{matrix} \qquad
  U = \begin{matrix}
    1&1&1&1&1&1&1&1\\
    2&3&3&3&3&3& & \\
    3&8&8& & & & & \\
    4& & & & & & & 
  \end{matrix}
\end{equation*}
\end{ex}

To prove Conjecture \ref{cat conj}, it is enough to show that there
is an involution $\theta$ on $\SSS$ such that
\begin{enumerate}
\item For every point $(w,T,U)\in\SSS$, we have
$\weight(\theta(w,T,U)) = \weight(w,T,U)$ (weight-preserving),
and if $(w,T,U)$ is not in the set $\SSS^\theta$ of 
fixed points of $\theta$, then
$\sign(\theta(w,T,U))=-\sign(w,T,U)$ (sign-reversing).
\item Each triple in the set $\SSS^\theta$ of fixed points of $\theta$ has
positive sign, and there is a bijection
$\SSS^\theta\rightarrow CT(\la;R)$ such that
$\weight(w,T,U) = q^{\charge(P)}$ where $(w,T,U)\mapsto P$.
\end{enumerate}

We define the cancelling involution by transforming the data
through a series of bijections that compute the
operator $\cat_{R_1}$, and then applying
a cancelling involution on the resulting set.
Denote the length of a word $u$ by $|u|$.

Given the triple $(w,T,U)\in\SSS'$, consider the sequence of words
given by the inverse image of the tableau pair $(T,U)$ under column RSK.
Here we use a nonstandard indexing for the weakly
increasing words, writing
\begin{equation*}
P(u^m u^{m-1} \dots u^1 u^n u^{n-1} \dots u^{m+1}) = T.
\end{equation*}
Since $U$ has content $(\beta(w),\alpha(w)-R_1)$,
this nonstandard indexing allows us to say that the length of
$u^i$ is the $i$-th part of the weight
$\xi(w)-(R_1,0^{n-m}) = (\alpha(w)-R_1,\beta(w))$.

\begin{ex} Continuing the above example, the words $u^i$ are given by
$u^2=56$, $u^1=u^8=u^7=\emptyset$, $u^6=8$,
$u^5=445688$, $u^4=4$, $u^3=33356777$.
\end{ex}

Now let
\begin{equation} \label{v def}
  v^i = \begin{cases}
  	i^{\gamma_i} u^i & \text{if $1\le i\le m$} \\
  	u^i	& \text{if $m<i\le n$.}
  	\end{cases}
\end{equation}
The word $v^i$ is weakly increasing since $u^i$ is
weakly increasing and consists of letters in the alphabet $[m+1,n]$.
Note that $|v^i| = \xi(w)_i$.

\begin{ex} In the above example we have $v^1=111$, $v^2=22256$, and
$v^i=u^i$ for $2<i\le n$.
\end{ex}

Finally, let $(P,Q)$ be the pair of tableaux given by the image under
column RSK, of the sequence of words $\{v^i\}$,
so that $P(v^n v^{n-1} \dots v^1) = P$.  Note that
the contents of $P$ and $Q$ are $\gamma$ and $\xi(w)$ respectively.

\begin{ex} The tableaux $P$ and $Q$ are given by
\begin{equation*}
P=\begin{matrix}
  1&1&1&5&5&6&7&7\\
  2&2&2&6&6&7& & \\
  3&3&3&8&8& & & \\
  4&4&4& & & & & \\
  5& & & & & & & \\
  8& & & & & & &
\end{matrix}  \qquad
Q=\begin{matrix}
1&1&1&2&2&3&3&3 \\
2&2&2&3&3&5& &  \\
3&3&3&5&5& & &  \\
4&5&5& & & & &  \\
5& & & & & & &  \\
6& & & & & & &
\end{matrix}
\end{equation*}
\end{ex}

Define the map $\Phi:\SSS'\rightarrow \Phi(\SSS')$ given by
$\Phi(w,T,U)=(w,P,Q)$.  
It is clear from the definitions that $\Phi$ is injective.
Define a sign and weight on $\Phi(\SSS')$
by $\sign(w,P,Q)=(-1)^w$ and $\weight(w,P,Q)=q^{\charge(P)}$.
By definition $\Phi$ is sign-preserving.
We now show that $\Phi$ preserves weight.

\begin{lem} \label{charge crank} In the above notation,
\begin{equation*}
  \charge(P) = \charge(T) + |\alpha(w)| - |R_1|
\end{equation*}
provided that $\gamma$ is a partition.
\end{lem}
\begin{proof} Using Theorem \ref{charge characterization} we
compute the charge of $P$.
\begin{equation*}
\begin{split}
  \charge(P) &= \charge(v^n v^{n-1} \dots v^1)  \\
    &= \charge(u^n u^{n-1} \dots u^{m+1}
    m^{\gamma_m} u^m (m-1)^{\gamma_{m-1}} u^{m-1} \dots 1^{\gamma_1} u^1) \\
    &= |u^1| + \charge(u^1 u^n u^{n-1} \dots u^{m+1} \\
     &\qquad
    m^{\gamma_m} u^m (m-1)^{\gamma_{m-1}} u^{m-1} \dots 2^{\gamma_2} u^2) \\
    &= |u^2| + |u^1| + \charge(u^2 u^1 u^n \dots u^{m+1}
    m^{\gamma_m} u^m \dots 3^{\gamma_3} u^3) \\
    &= \vdots \\
    &= |u^m|+\dots+|u^1| + \charge(u^m \dots u^1 u^n \dots u^{m+1}) \\
    &= |\alpha(w)|-|R_1| + \charge(T).
\end{split}
\end{equation*}
\end{proof}

We define a sign-reversing, weight-preserving involution
$\theta'$ on the set $\Phi(\SSS')$.  Let $(w,P,Q)\in\Phi(\SSS')$.
There are two cases.
\begin{enumerate}
\item $Q$ is lattice.  It follows that $Q$ has partition content
and $w$ is the identity.  Define $\theta'(w,P,Q)=(w,P,Q)$.
\item $Q$ is not lattice.  Let $r+1$ be the rightmost
letter in the row-reading word of $Q$ that violates the lattice
condition.  Define $\theta'(w,P,Q)=(w',P',Q')$ where
$w' = w s_r$, $P'=P$, and $Q' = s_r e_r Q$.
\end{enumerate}

\begin{ex} In computing $\theta'(w,P,Q)=(w',P',Q')$, the first
violation of latticeness in the row-reading word of
$Q$ occurs at the cell $(1,8)$ so $r=2$.  We have
$w'=31254678$, $P'=P$ and
\begin{equation*}
Q'=\begin{matrix}
1&1&1&2&2&2&2&3 \\
2&2&2&3&3&5& &  \\
3&3&3&5&5& & &  \\
4&5&5& & & & &  \\
5& & & & & & &  \\
6& & & & & & &
\end{matrix}
\end{equation*}
\end{ex}

\begin{lem} \label{inv lem}
$\theta'$ is a sign-reversing, weight-preserving involution
on the set $\Phi(\SSS')$.
\end{lem}
\begin{proof} By definition $\theta'$ is sign-reversing and weight-
preserving.  By Lemma \ref{lattice} $\theta'$ is an involution.
It remains to show that $\theta'$ stabilizes the set $\Phi(\SSS')$.
Let $(w,P,Q)\in\Phi(\SSS')$ and $\theta'(w,P,Q)=(w',P',Q')$.
We may assume $(w,P,Q)$ is not a fixed point of $\theta'$.
Let $v'$ be to $(P',Q')$ as $v$ is to $(P,Q)$ in the definition of
$\Phi$.  It is enough to show that $(v')^{i}$ starts with the subword
$i^{\gamma_i}$ for every $1\le i\le m$, since the other steps in the map
$\Phi$ are invertible by definition.

Since $Q$ and $Q' = s_r e_r(Q)$ are in the same $r$-string,
Lemma \ref{two row dual} applies.  
There is nothing to prove unless $r\le m$.
Suppose first that $r<m$.  We need only check that
${v'}^r$ starts with $r^{\gamma_r}$ and
${v'}^{r+1}$ starts with $(r+1)^{\gamma_{r+1}}$.
Since $v^r = r^{\gamma_r} u^r$ and
$v^{r+1} = (r+1)^{\gamma_{r+1}} u^{r+1}$ where
all the letters of $u^r$ and $u^{r+1}$ are strictly greater than $m$,
it follows that
\begin{equation*}
  P({v'}^{r+1} {v'}^r)|_{[r,r+1]} =
  P(v^{r+1} v^r)|_{[r,r+1]} = P((r+1)^{\gamma_{r+1}} r^{\gamma_r})
\end{equation*}
This, together with the fact that ${v'}^{r+1}$ and ${v'}^r$ are weakly
increasing words, implies that all of the letters $r+1$ must precede all
of the letters $r$ in the word ${v'}^{r+1} {v'}^r$, that is,
${v'}^i$ starts with $i^{\gamma_i}$ for $i\in[r,r+1]$.

The remaining case is $r=m$.  Then $v^r=r^{\gamma_r} u^r$
and $v^{r+1} = u^{r+1}$.  Let us calculate ${v'}^{r+1}$ and ${v'}^r$
using a two-row jeu-de-taquin.  Let $R$ (resp. $R'$) be the (skew) two row
tableau with first row $v^r$ (resp. ${v'}^r$ and second row $v^{r+1}$
(resp. ${v'}^{r+1}$ in which the two rows achieve the maximum overlap.
The overlaps of $R$ and $R'$ are equal by Lemma \ref{overlap} and
the fact that $Q$ and $Q'$ are in the same $r$-string
and hence have the same $r$-paired letters.  Furthermore this
common overlap is at least $\gamma_r$.  To see this, note that
the overlap weakly exceeds the minimum of $\gamma_r$ and $|u^{r+1}|$
since all of the letters in $u^{r+1}$ have values in the alphabet
$[m+1,n]=[r+1,n]$ and there are $\gamma_r$ copies of $r$ in $u^r$.
On the other hand, $|u^{r+1}|>\gamma_r$, for otherwise by
Lemma \ref{overlap} all of the letters $r+1$ in $Q$ would be
$r$-paired, contradicting the choice of $r$.

We calculate $R'$ from $R$ in two stages.
Let $R"$ be the two row skew tableau (whose rows have maximum
overlap) such that $P(R")=P(R)$, where the first row of $R"$ is
one cell longer than that of $R$.  By Lemma \ref{two row dual} this
tableau exists since $Q$ has an $r$-unpaired letter $r+1$; the
corresponding recording tableau is $e_r Q$.  Furthermore $R"$ is
obtained by sliding the ``hole" in the cell just to the left of the first
letter in the first row of $R$, into the second row.
By the same reasoning as above, $R"$ has the same overlap that $R$ does.
Finally we calculate $R'$ from $R"$ by another two row jeu-de-taquin.
If the first row of $R"$ is shorter than the second, we are done,
for in this case the first row ${v'}^r$ of $R'$ contains the first
row of $R"$, which in turn contains the first row of $R$, which
contains $r^{\gamma_r}$.  So suppose
the second row of $R"$ is shorter than the first, by $p$ cells, say.
Now $p$ is less than or equal to the number of cells on the right
end of the first column of $R"$ that have no cell of $R"$ below them.
Since $R"$ has maximum overlap it follows that
when $p$ holes are slid from the second row of $R"$ to the first,
they all exchange with numbers lying in the portion of the first row
of $R"$ that extends properly to the right of the second.
Thus the subword $r^{\gamma_r}$ remains in the first row of $R'$,
and we are done.
\end{proof}

\begin{ex} In $v'$ all the subwords are the same as in $v$ except
that ${v'}^2 = 2225677$ and ${v'}^3=333567$.  In this example $r=m$.
The tableaux $R$, $R"$, and $R'$ are given below.
\begin{equation*}
\begin{split}
R&=\begin{matrix}
  \sq&\sq&\sq&2&2&2&5&6 \\
  3&3&3&5&6&7&7&7
\end{matrix} \\
R"&=\begin{matrix}
  \sq&\sq&2&2&2&5&6&7 \\
  3&3&3&5&6&7&7&
\end{matrix} \\
R'&=\begin{matrix}
  \sq&2&2&2&5&6&7&7 \\
  3&3&3&5&6&7& &
\end{matrix}
\end{split}
\end{equation*}
\end{ex}

Thus we may define a sign-reversing, weight-preserving involution
$\theta$ on $\SSS'$ by $\theta = \Phi^{-1} \circ \theta' \circ \Phi$.
By definition the fixed points of $\theta'$ are the triples
$(w,P,Q)$ where $Q$ is the unique column strict tableau of shape and
content $\la$, $w$ is the identity, and $P$ is a column strict tableau
of content $\gamma$ such that $P|_{A_1} = Y_1$.

\begin{lem} \label{reduction}
Conjecture \ref{cat conj} holds provided that the
involution $\theta'$ stabilizes the subset $\Phi(\SSS)$
of $\Phi(\SSS')$.
\end{lem}
\begin{proof}  Suppose that $\theta'$ stabilizes $\Phi(\SSS)$.
Equivalently, the involution $\theta$ stabilizes $\SSS$.
Since $\Phi$ is sign and weight-preserving, the generating function
of the fixed points $\SSS^\theta$
of the restriction of $\theta$ to $\SSS$
is the same as that of the set $\Phi(\SSS)^{\theta'}$
But it is easy to see that $\SSS^\theta$
is precisely the triples
$(id,P,Q)$ where $Q$ is the unique column strict tableau of shape and
content $\la$ and $P$ is $R$-catabolizable.
\end{proof}

In fact, the proof of Conjecture \ref{cat conj} reduces to
the case of a non-fixed point $(w,T,U)$ with $r=m$.

\begin{lem} In the notation of Lemma \ref{reduction} and the
definition of $\theta$, if $r\not=m$ then $T'=T$.
\end{lem}
\begin{proof} Let $u'$ be to $(T',U')$ as $u$ is to $(T,U)$ in the
definition of $\Phi$.  It follows from Lemma \ref{two row dual} that
${u'}^i =u^i$ for $i\not\in[r,r+1]$ and
$P({u'}^{r+1}{u'}^r)=P(u^{r+1} u^r)$.  Since $r\not=m$ we have
\begin{equation*}
\begin{split}
  P({u'}^m \dots {u'}^1) &= P(u^m \dots u^1) \\
  P({u'}^n \dots {u'}^{m+1}) &= P(u^n \dots u^{m+1}) \\
\end{split}
\end{equation*}
which implies that
\begin{equation*}
\begin{split}
  T' &= P({u'}^m \dots {u'}^1 {u'}^n \dots {u'}^{m+1}) \\
 &= P(u^m \dots u^1 u^n \dots u^{m+1}) \\
 &= T
\end{split}
\end{equation*}
\end{proof}

We believe that $T'$ is always $\Rhat$-catabolizable,
but can only prove it in the following
cases, where the definition of catabolizability becomes
quite simple.  The following result is an immediate
consequence of the definitions.

\begin{prop} \label{hook cat} Let $\eta$ be a hook partition
(that is, $\eta_i=1$ for $i>1$).
Then $S$ is $R$-catabolizable if and only if $S|_{A_1} = Y_1$
and $S$ has content $\gamma$.
\end{prop}

\begin{cor} Conjecture \ref{cat conj} holds when $\eta=(m,1^{n-m})$.
\end{cor}
\begin{proof} By Proposition \ref{hook cat} the tableau
$T$ is $\Rhat$-catabolizable in the alphabet $[m+1,n]$
if and only if $T$ has content $(0^m,\gammahat)$.  In this case
the sets $\SSS$ and $\SSS'$ coincide.
\end{proof}

\begin{cor} \label{two part cor}
Conjecture \ref{cat conj} holds when $\eta$ has two parts.
\end{cor}
\begin{proof} Consider the first recurrence \eqref{recurrence}, which
by \eqref{base of recurrence} takes the form
\begin{equation} \label{two part rec}
\begin{split}
  \K_{\la;(R_1,R_2)}(q) &= \sum_{w\in W/(W_y\times W_z)}
  (-1)^w q^{|\alpha(w)|-|R_1|} \sum_\sigma
	\LR^\sigma_{\alpha(w)/R_1,\beta(w)} \delta_{\sigma,R_2} \\
	&= \sum_{w\in W/(W_y\times W_z)} (-1)^w q^{|\alpha(w)|-|R_1|} 
	\LR^{R_2}_{\alpha(w)/R_1,\beta(w)}
\end{split}
\end{equation}
Suppose that the $w$-th term is nonzero and $w$ is not the identity.
Necessarily $\alpha(w) \supseteq R_1$.  By examining the Jacobi-Trudi
determinant for $s_{\la/R_1}$, it must be the case that
$\la_{m+1}>\gamma_m$.  On the other hand we have
$\beta(w)_1>\beta(id)_1 = \la_{m+1}$.  Putting the inequalities together,
$\beta(w)_1>\gamma_m\ge\gamma_{m+1}$ by the dominance of the weight
$\gamma$.  But $\gamma_{m+1}$ is the first part of the partition $R_2$.
This means that $R_2$ cannot contain $\beta(w)$, so that the
LR coefficient vanishes, contradicting our assumption on $w$.

Equation \eqref{two part rec} now takes the form
\begin{equation} \label{two part Poincare}
  \K_{\la;(R_1,R_2)}(q) = q^{|\alpha|-|R_1|} \LR^{R_2}_{\alpha/R_1,\beta}
\end{equation}
where $\alpha=\alpha(id)$ and $\beta=\beta(id)$ are the first $m$
and last $n-m$ parts of $\la$ respectively.  Note that the skew
shape $\la/R_1$ consists of two disconnected skew shapes, namely
$\alpha/R_1$ and $\beta$.  Let $P$ be an $R$-catabolizable tableau of shape
$\la$ and $Q$ the unique column strict tableau of shape and content $\la$.
Let $T$ and $U$ be the tableaux such that $\Phi(id,T,U) = (id,P,Q)$.
For $(id,P,Q)$ to be in the image of $\Phi$, it must be shown that
$v^i$ starts with the subword $i^{\gamma_i}$ for $1\le i\le m$; but
this holds since $v^i$ is the $i$-th row of the tableau $P$,
and $P|_{A_1} = Y_1$ by the $R$-catabolizability of $P$.
By Lemma \ref{charge crank} $\charge(P) = |\alpha| - |R_1|$, since
$T=Y_2$ by the $R$-catabolizability of $P$ and the fact that
Yamanouchi tableaux have zero charge.  It follows from
Theorem \ref{fitting} that the map $P\mapsto U$ is a bijection
from $CT(\la;R)$ to a set of tableaux of cardinality
$\LR^{R_2}_{\alpha/R_1,\beta}$.  The result follows.
\end{proof}

\section{An application of catabolizability}
\label{app sec}

We prove a beautiful formula for the cocharge Kostka-Foulkes
polynomials stated by Lascoux \cite{La}.  The proof also settles Conjecture
\ref{cat conj} in the case that $\gamma=(1^t)$ and $\eta$ is a partition,
that is, when all rectangles are single columns of weakly decreasing
heights.  This entails a deeper study of catabolizability.

\subsection{Lascoux' formula}
For the pair of partitions $\la$ and $\mu$,
the cocharge Kostka-Foulkes polynomial $\tK_{\la,\mu}(q)$ is defined by
\begin{equation*}
  \tK_{\la,\mu}(q) = q^{n(\mu)} K_{\la,\mu}(q)
\end{equation*}
where $n(\mu) = \sum_i (i-1)\mu_i$.  For a combinatorial
description, define the cocharge of a word $u$ of partition
content $\mu$ by
\begin{equation*}
 \cocharge(u) = n(\mu) - \charge(u).
\end{equation*}
Then Theorem \ref{charge thm} can be rephrased as
\begin{equation} \label{CST cocharge}
  \tK_{\la,\mu}(q) = \sum_T q^{\cocharge(T)}
\end{equation}
where $T$ runs over the set of column strict tableaux
of shape $\la$ and content $\mu$.

These polynomials satisfy the monotonicity property
\begin{equation} \label{monotonicity}
\tK_{\la,\mu}(q) \le \tK_{\la,\nu}(q)
\end{equation}
coefficientwise provided that $\mu\dom\nu$.

To exhibit a combinatorial proof of the inequality
\eqref{monotonicity}, Lascoux associates to each standard tableau $S$
a partition $\cattype(S)$ and asserts that
\begin{equation} \label{standard cocharge}
  \tK_{\la,\mu}(q) = \sum_S q^{\cocharge(S)}
\end{equation}
where $S$ runs over the set of standard tableaux
of shape $\la$ and $\cattype(S) \dom \mu$.

Define the tableau operator $\Cat(S)=H_1(S)$, which may be computed by
column-inserting the first row of $S$ into the remainder of $S$.
Let $d_1(S)$ be the maximum number $i$ such that the numbers $1$
through $i$ are all in the first row of $S$.  Then $\cattype(S)$
is the sequence of intergers whose first part is $d_1(S)$ and
whose $i$-th part is given by $d_1(\Cat^i(S))-d_1(\Cat^{i-1}(S))$.
It can be shown that $\cattype(S)$ is a partition.

\begin{ex} The powers of $\Cat$ on the standard tableau $S$
are given below.
\begin{equation*}
\begin{split}
S &= \begin{matrix}
1&2&3&4&7\\
5&6&9& & \\
8& & & & 
\end{matrix}
\qquad
\Cat(S) = \begin{matrix}
1&2&3&4&5&6&9\\
7&8& & & & & \\
 & & & & & & \\
\end{matrix} \\
\Cat^2(S) &= \begin{matrix}
1&2&3&4&5&6&7&8\\
9& & & & & & & \\
\end{matrix}
\qquad
\Cat^3(S) = \begin{matrix}
1&2&3&4&5&6&7&8&9\\
 & & & & & & & &
\end{matrix}
\end{split}
\end{equation*}
So the sequence of $d_1$ of the powers of $\Cat$ on $S$ is
$(4,6,8,9)$ and $\cattype(S)=(4,2,2,1)$.
\end{ex}

Let $\alpha$ be a sequence of nonnegative integers and
$\TTT(\alpha)$ the set of column strict tableaux
of content $\alpha$ and arbitrary partition shape.
To prove Lascoux' formula, in light of
\eqref{CST cocharge} and \eqref{standard cocharge}
it clearly suffices to exhibit an embedding
$\theta_\mu:\TTT(\mu)\rightarrow \TTT((1^n))$ (where $n=|\mu|$)
that is shape- and cocharge-preserving, and satisfies the
additional property that $S$ is in the image of $\theta_\mu$
if and only if $\cattype(S)\dom\mu$.  Lascoux defined such an
embedding but did not give a proof of the characterization
of the image of $\theta_\mu$.  We supply a proof of this
last fact.

\subsection{Cyclage and canonical embeddings}

The embeddings $\theta_\mu$ are best understood in terms
of the cyclage poset structure \cite{La} \cite{LS2}.
The \textit{cyclage} is the covering relation of a graded
poset structure on $\TTT(\alpha)$.  It is defined as follows.

Suppose first that $\alpha$ is a partition $\mu$.
For two tableaux $T$ and $S$ in $\TTT(\mu)$, say that
$T$ covers $S$ if there exists a letter $a>1$ and
a word $x$ such that $P(ax)=T$ and $P(xa)=S$.
Equivalently, $T$ covers $S$ if there is a corner cell
$s$ of $T$ such that the reverse column insertion on $T$
at $s$ results in a letter $a>1$ and a tableau $U$,
and the row insertion of the letter $a$ into $U$ produces $S$.

The above relation (called the cyclage)
is the covering relation of a partial order since
$\cocharge(T) = \cocharge(S) + 1$ when $T$ covers $S$,
by Theorem \ref{charge characterization}.

\begin{ex} Starting with the tableau $T$ and an underlined cell $s$,
the pair $(x,U)$ and the tableau $S$ are computed.
\begin{equation*}
\begin{split}
T &= \begin{matrix}
  1&1&1&2&3\\
  2&3&\underline{4}& & \\
  4& & & & 
\end{matrix} \qquad
x=2\qquad
U=\begin{matrix}
  1&1&1&2&3\\
  3&4& & & \\
  4& & & & 
\end{matrix} \\
S&=\begin{matrix}
  1&1&1&2&2\\
  3&3& & & \\
  4&4& & & 
\end{matrix}
\end{split}
\end{equation*}
\end{ex}

Now let us return to the case of arbitrary content $\alpha$.
Let $\mu = \alpha^+$ and $w=(w_\alpha)^{-1}$ the permutation
such that $w \alpha = \mu$.  Then $w$ defines a bijection
$\TTT(\alpha)\rightarrow\TTT(\mu)$ given by
$T\mapsto w T$ (see \eqref{plactic action}) that is shape-preserving.

Say that $T \ge S$ for $T,S\in\TTT(\alpha)$ if
$wT \ge wS$, that is, the partial order on $\TTT(\alpha)$ is defined
to make $w$ into a poset isomorphism.

\begin{thm} \cite{LS2} The cyclage endows $\TTT(\alpha)$ with the
structure of a graded poset with grading given by cocharge.
The unique bottom element of $\TTT(\alpha)$ is the one-row tableau
of content $\alpha$.
\end{thm}

The family of posets $\{\TTT(\alpha)\}$ for varying alpha,
is equipped with functorial grade-preserving embeddings.
For two compositions $\alpha$ and $\beta$, say that
$\alpha \dom\beta$ if $\alpha^+ \dom \beta^+$.
Suppose $\alpha\dom\beta$.  There is an embedding of posets
$\theta^{\beta}_\alpha:\TTT(\alpha)\rightarrow\TTT(\beta)$
that can be defined as follows.

First, if $\alpha^+ = \beta^+$, choose any $w\in W$ such that
$w \alpha =\beta$ and let $\theta^\beta_\alpha = w$.

Second, if $\beta_i = \alpha_i$ for all $i>2$,
$\beta_1 = \alpha_1-1$ and $\beta_2 = \alpha_2 + 1$
where $\alpha_1>\alpha_2 + 1$, then let
$\theta^{\beta}_{\alpha} = f_1$, the crystal lowering operator,
which in this case merely changes the rightmost letter $1$ in a
tableau to a $2$.

Now let $\alpha \dom \beta$.  Then there is a sequence
$\alpha=\gamma^0,\gamma^1,\dots,\gamma^p=\beta$ of 
compositions, where either $(\gamma^i)^+ = (\gamma^{i+1})^+$
or $\gamma^{i+1} = \gamma^i + (-1,1,0,\dots)$.  Define
\begin{equation*}
  \theta^\beta_\alpha = \theta^{\gamma^p}_{\gamma^{p-1}}
  	\circ \dots \circ \theta^{\gamma^1}_{\gamma^0}.
\end{equation*}

\begin{thm} \label{embeddings} \cite{La} Suppose $\alpha\dom\beta$.
The map $\theta^\beta_\alpha$ is independent of the
sequence of compositions $\{\gamma^i\}$ and is
shape-preserving and an embedding of graded posets.  Furthermore if
$\beta\dom\gamma$ then
$\theta^{\gamma}_\alpha = \theta^\gamma_\beta \circ \theta^\beta_\alpha$.
\end{thm}

For $|\alpha| = n$, denote by $\theta_\alpha = \theta^{(1^n)}_\alpha$
the embedding of $\TTT(\alpha)$ into the standard tableaux $\TTT((1^n))$.
Lascoux gives the following characterization of the image
of $\theta_\mu$.

\begin{thm} \label{cyc image}
Let $\mu$ be a partition of $n$.  Then
the image of $\theta_\mu$ is the set of standard tableaux
$S$ such that $\cattype(S) \dom\mu$.
\end{thm}

\subsection{Proof of Theorem \ref{cyc image}}
\label{cyc std proof}

Our proof uses several reformulations of the condition
$\cattype(S)\dom\mu$.
Let $Z_m$ be the one-row standard tableau given by the numbers
from $1$ to $m$.  Suppose the standard tableau $S$ contains $Z_m$.
Define $\Cat_m(S) = H_1(S-Z_m)$.
Write $\mu=(m,\muhat)$.  The following definition is not
consistent with the previous notion of $R$-catabolizability
for any $R$ since it slices the tableau in the wrong direction.
Say that $S$ is \textit{$\mu$-catabolizable} if
$S$ and $\mu$ are both empty, or if $S$ contains
$Z_m$ and $\Cat_m(S)$ is $\muhat$-catabolizable
in the alphabet $[m+1,n]$.

\begin{prop} \label{cat and type}
$S$ is $\mu$-catabolizable if and only if
$\cattype(S)\dom\mu$.
\end{prop}
\begin{proof} This is trivial if $\mu$ has one part.
So let $\mu=(m,\muhat)$.  We may assume that $S$ contains
$Z_m$ for otherwise both conditions on $S$ fail.
From the definitions one has
$\Cat(S) = P(Z_m \Cat_m(S))$.  In particular, since all the
letters of $Z_m$ are smaller than all of those in $\Cat_m(S)$,
one obtains $\Cat(S)$ from $\Cat_m(S)$ by pushing the first row
to the right by $m$ cells and placing $Z_m$ in the vacated positions.
In view of this, it is clear that 
$d_1(\Cat^i(\Cat_m(S)))=d_1(\Cat^{i+1}(S))$ for all $i$.
So if we write $\cattype(S) = \nu$ then
$\cattype(\Cat_m(S)) = (\nu_1-\mu_1+\nu_2,\nu_3,\nu_4,\dots)$.
By induction $\Cat_m(S)$ is $\muhat$-catabolizable if and only if
$\cattype(\Cat_m(S)) \dom \muhat$,
which is equivalent to $\nu\dom\mu$.
\end{proof}

The next reformulation of the condition $\cattype(S)\dom\mu$
is precisely the transpose of the condition of
$R$-catabolizability where $R_i=(1^{\mu_i})$.

Define the vertical slicing operator $V_c$ as follows.
For the (skew) tableau $T$, let $V_c(T)=P(T_e T_w)$
where $T_e$ and $T_w$ are the east and west subtableaux
obtained by slicing $T$ vertically between the $c$-th and
$(c+1)$-st columns.

Suppose that $S$ contains $Z_m$.  Define the operator
$\CCat_m(S) = V_m(S-Z_m)$.
Say that $S$ is \textit{$\mu$-column catabolizable} if
$S$ and $\mu$ are empty, or if $S$ contains $Z_{\mu_1}$
and $\CCat_m(S)$ is $\muhat$-column catabolizable.

The following equivalence is far from obvious.

\begin{prop} \label{row and col cat}
A standard tableau is $\mu$-catabolizable
if and only if it is $\mu$-column catabolizable.
\end{prop}

To prove this it is necessary to
consider the following restrictions of the
cyclage relation for standard tableaux.  Let us use
the notation $T \ge S$ for the cyclage partial order.
Say that $T$ covers $S$ in the partial order
$\ge_{(r,)}$ if, $T$ covers $S$ in the order $\ge$ and
(in the notation of the definition of $\ge$) the ``starting cell''
$s$ lies in a row strictly below the $r$-th.

Say that $T$ covers $S$ in the order $\ge_{(,c)}$ if
$T$ covers $S$ in the order $\ge$ and the ``ending cell''
$s'$ (defined by the difference of the shapes of $S$ and
the intermediate tableau $U$) is in a column strictly to the right
of the $c$-th.  Another viewpoint is to reverse the cyclage construction
(call this cocyclage).  One starts with a cell $s'$, performs a reverse
row insertion on $S$ at $s'$ to produce a letter $x$ and a tableau $U$,
then column inserts $x$ into $U$ to produce the tableau $T$.  
Then if $T \ge_{(,c)} S$, the starting cell $s'$ of the cocyclage
must start at a cell in a column strictly right of the $c$-th.

These orders are compatible with the two notions of catabolizability.

\begin{lem} \label{row cat and cyc}
If $T$ is $\mu$-catabolizable and $T \ge_{(1,)} S$,
then $S$ is also.
\end{lem}
\begin{proof} Without loss assume that $T \ge_{(1,)} S$ is a covering
relation.  Let $\mu=(m,\muhat)$.  Let $s$ be the cell where the cyclage
starts, $U$ and $a$ as in the definition of cyclage.  By definition we have
$T = P(a U)$ and $S=P(Ua)$.
Let $\Th$ be obtained by removing the first row from $T$.  Define
$\Sh$ and $\Uh$ similarly.  
Since $T$ is $\mu$-catabolizable, the first row of $T$ has the form
$Z_m x$ where $x$ is a word.  Since $s$ is not in the first row and
the bumping paths of reverse column insertions move weakly south,
it follows that the first rows of $T$ and $U$ coincide and
$a \Uh \Kn \Th$.  There are two cases:
\begin{enumerate}
\item $xa$ is a weakly increasing word.  In this case the first row of $S$
is equal to $Z_m x a$ and $\Sh=\Uh$, so that
\begin{equation*}
H_1(S-Z_m) \Kn x a \Sh = x a \Uh \Kn x U \Kn H_1(T-Z_m).
\end{equation*}
But $H_1(T-Z_m)$ is $\muhat$-catabolizable by definition.
\item $xa$ is not weakly increasing.
Here the letter $a$ is strictly smaller
than some letter of $x$.  Let $y$ the weakly increasing word and $b$
the letter such that $xa \Kn by$.  Then $Z_m y$ is the first row
of $S$ and $\Sh \Kn \Uh b$.  We have
\begin{equation*}
\begin{split}
H_1(T-Z_m) &\Kn x \Th \Kn x a \Uh \Kn b y \Uh \\
H_1(S-Z_m) &\Kn y \Sh \Kn y \Uh b
\end{split}
\end{equation*}
By induction it is enough to show that
$P(b y \Uh) \ge_{(1,)} P(y \Uh b)$.

Let $s"$ be the cell where the column insertion of
$b$ into $P(y \Uh)$ ends.  It is enough to show that $s"$ is not
in the first row.  Let $y= cz$, where $c$ is the first letter of $y$
and $z$ the remainder of $y$.  Since $xa$ is not weakly increasing and
and $xa \Kn by$, we have $b>c$.  Thus the
bumping path of the column insertion of $b$ into
$P(cz \Uh)$ is strictly south of that of the column insertion of $c$ into
$P(z \Uh)$.  It follows that $s"$ is not in the first row.
\end{enumerate}
\end{proof}

\begin{lem} \label{col cat and cyc}
If $S$ is $\mu$-column catabolizable and
$T \ge_{(,m)} S$, then $T$ is also.
\end{lem}
\begin{proof} Assume that $T \ge_{(,m)} S$ is a covering relation.
Let $s'$ be the starting cell of the cocyclage from $S$ to $T$,
and $a$ and $U$ as in the definition of cyclage.  We have
$S=P(Ua)$ and $T=P(aU)$.  Let $\mu=(m,\muhat)$.
Let $S_e$ and $S_w$ be the east and west subtableaux obtained by
slicing $S$ vertically just after its $m$-th column.
Define $T_e$, $T_w$, $U_e$, and $U_w$ similarly.
Then $S_e \Kn U_e a$.  Now the bumping path of a reverse row insertion
moves weakly east, so $S_w=U_w$.  There are two cases.
\begin{enumerate}
\item The column insertion of $a$ into $U$ ends weakly west of the
$m$-th column.  In this case $T_w \Kn a U_w$ and $T_e = U_e$.
It follows that
\begin{equation*}
V_m(T-K) \Kn T_e T_w \Kn U_e a U_w \Kn S_e S_w \Kn V_m(S-K).
\end{equation*}
Thus $\CCat_m(T)=\CCat_m(S)$ and $T$ is $\mu$-column catabolizable.
\item The column insertion of $a$ into $U$ ends in a column strictly
east of the $\mu_1$-th.  Let $b$ be the letter that is bumped from
the $\mu_1$-th column to the $\mu_1+1$-st during this column insertion.
Then $a U_w \Kn T_w b$ and $b U_e\Kn T_e$.  We have
\begin{equation*}
\begin{split}
V_m(T-K) &\Kn T_e T_w \Kn b U_e T_w \\
V_m(S-K) &\Kn S_e S_w \Kn U_e a U_w \Kn U_e T_w b
\end{split}
\end{equation*}
It is enough to show that $P(b U_e T_w) \ge_{(,m)} P(U_e T_w b)$.
Let $s"$ be the cell where the row insertion of $b$ into $P(U_e T_w)$
ends.  It is enough to show that $s"$ lies strictly east of the
$\mu_1$-th column.  By the assumption of this case,
$T_w b$ is a tableau of partition shape.
By \cite{White} it follows that
$s"$ lies strictly to the east of all of the cells of the skew shape
given by $\shape(P(U_e T_w))/\shape(U_e)$.  But this skew shape
contains a horizontal strip of length $m$ (since the first
row of $T_w$ has length $m$), so it follows that $s"$
is strictly east of the $m$-th column.
\end{enumerate}
\end{proof}

\begin{lem} Suppose $T$ contains $Z_m$.  Then
\begin{equation*}
\begin{split}
  \Cat_m(T) &\le_{(1,)} \CCat_m(T) \\
  \Cat_m(T) &\le_{(,m)} \CCat_m(T)
\end{split}
\end{equation*}
\end{lem}
\begin{proof}  Slice $T$ horizontally between the first and second rows
and vertically between the $m$-th and $m+1$-st columns.
Let the northwest, northeast, southwest, and southeast subtableaux be
denoted by $Z_m$, $T_{ne}$, $T_{sw}$, and $T_{se}$ respectively.
By the definitions, we have
\begin{equation*}
\begin{split}
  \Cat_m(T)&=H_1(T-Z_m) \Kn T_{ne} T_{sw} T_{se} \\
  \CCat_m(T) &= V_m(T-Z_m) \Kn T_{se} T_{ne} T_{sw}.
\end{split}
\end{equation*}
It is enough to show that the relation
\begin{equation*}
  P(T_{ne} T_{sw} T_{se}) \le P(T_{se} T_{ne} T_{sw})
\end{equation*}
(which holds since all of the numbers in the tableau $T_{se}$
are strictly greater than $2m$), also holds in
the orders $\le_{(,m)}$ and $\le_{(1,)}$.

If the tableau $T_{se}$ is column inserted into the tableau $T_{ne}$,
all of the bumping paths end in rows strictly south of the first.
It follows that if one first row inserts $T_{sw}$ into $T_{ne}$,
and then column inserts tableau $T_{se}$, the column insertions
are pushed weakly to the south and hence must still end in rows strictly
south of the first.  That means there is a sequence of cocyclages
starting in rows strictly south of the first, that prove the relation
\begin{equation*}
P(T_{se} T_{ne} T_{sw}) \ge_{(1,)} P(T_{ne} T_{sw} T_{se})
\end{equation*}
The proof for $\ge_{(,m)}$ is similar.
\end{proof}

Finally the proof of Proposition \ref{row and col cat} is given.
\begin{proof} Each of the statements implies the next,
using induction and the above lemmas.
\begin{enumerate}
\item $T$ is $\mu$-column catabolizable.
\item $T$ contains $Z_m$ and $\CCat_m(T)$ is $\muhat$-column catabolizable.
\item $T$ contains $Z_m$ and $\CCat_m(T)$ is $\muhat$-catabolizable.
\item $T$ contains $Z_m$ and $\Cat_m(T)$ is $\muhat$-catabolizable.
\item $T$ is $\mu$-catabolizable.
\end{enumerate}
Conversely, replace (3) with the following:
$T$ contains $Z_m$ and $\Cat_m(T)$ is $\muhat$-column catabolizable.
Then each statement implies the previous one.
\end{proof}

For the rest of the proof of Theorem \ref{cyc image},
a few more lemmas are needed.

\begin{lem} \label{long first row}
Suppose $S'\in\TTT(\mu)$ with first row
$1^m$ and remainder $\Sh'$ (of content
$(0,\mu_2,\mu_3,\dots)$).  Let $S=\theta_\mu(S')$.
Then the first row of $S$ is the tableau $Z_m$ and the rest
is given by $\Sh$, where $\Sh = \theta'_{\muhat}(\Sh')$
and $\theta'_{\muhat} = \theta_{(0,\mu_2,\mu_3,\dots)}^{(0^m,1^{n-m})}$.
\end{lem}
\begin{proof} For a tableau $T$ let $k+T$ denote the tableau
obtained by adding the integer $k$ to every letter in $T$.
The lemma is proven by careful explicit calculation
of the map $\theta_\mu$ and $\theta'_{\muhat}$.
We compute the map $\theta_\mu$ on the tableau $S'$
by the composition of the maps that change content as follows:
\begin{equation*}
  \mu \rightarrow (\muhat,0^{n-m-1},m) \rightarrow (1^{n-m},m)
  \rightarrow (m,0^{m-1},1^{n-m}) \rightarrow (1^m,1^{n-m}).
\end{equation*}
The first and third are given by crystal permutation operators
and the second and fourth by $\theta_{\muhat}$ and $\theta_{(m)}$
respectively.  By direct computation the image of $S'$ under the first map
is obtained by placing the number $(n-m+1)$ at the bottom of each of the
first $m$ columns of the tableau $-1+\Sh'$.  Then the second map
leaves the letters $(n-m+1)$ alone and applies $\theta_{\muhat}$
to the letters of the subtableau $-1+\Sh'$.  Again by direct computation
the third map produces the tableau with first row $1^m$ and remainder
$m+\theta_{\muhat}(-1+\Sh')$.  The fourth map leaves this remainder
alone and changes the first row from $1^m$ to $Z_m$.
In particular $\Sh=m+\theta_{\muhat}(-1+\Sh')$.
Now the map $\theta'_{\muhat}$ on the tableau $\Sh'$ can be computed
using the composition of the maps that change content as follows:
\begin{equation*}
  (0,\muhat) \rightarrow \muhat \rightarrow (1^{n-m})
  \rightarrow (0^m,1^{n-m}).
\end{equation*}
Now the first map is the crystal permutation operator, that
produces the tableau $-1+\Sh'$.  This kind of trivial relabelling occurs
when a crystal reflection operator $s_r$ acts on a word
that has no $r$'s or no $(r+1)$'s.  The second map is $\theta_{\muhat}$,
which produces $\theta_{\muhat}(-1+\Sh')$.  The third map is another
crystal permutation operator that produces $m+\theta_{\muhat}(-1+\Sh')$.
We have shown that
\begin{equation*}
  \Sh=m+\theta_{\muhat}(-1+\Sh')=\theta'_{\muhat}(\Sh').
\end{equation*}
\end{proof}

\begin{lem} \label{theta and north hat}
Suppose $X'\in\TTT(0,\muhat)$,
$X\in\TTT(0^m,1^{n-m})$, $S'=P(X' 1^m)$, and
$S=P(X Z_m)$.  Then in the notation of the previous lemma,
$\theta'_{\muhat}(X')=X$ if and only if
$\theta_\mu(S')=S$.
\end{lem}
\begin{proof} The proof proceeds by induction on charge.
Let $X_e$ and $X_w$ be the east and west subtableaux obtained by
slicing $X$ vertically between the $m$-th and $(m+1)$-st columns.
Make similar notation for the tableaux $X'$, $S$, and $S'$.
Since the letters of the word $1^m$ (resp. $Z_m$) are strictly
smaller than those in $X'$ (resp. $X$), it follows that
$S'$ (resp. $S$) is obtained from $X'$ (resp. $X$)
by pushing the subtableau $S'_w$ (resp. $S_w$)
down by a row, placing the word $1^m$ (resp. $Z_m$) in the
vacated cells, and leaving the other subtableau $S'_e$
(resp. $S_e$) in place.

So without loss we may assume that $X'$ and $X$ have the same shape
(and hence $S'$ and $S$ do as well).
Suppose first that $X'$ has at most $m$ columns.  Then
$X'=\Sh'$ and $X=\Sh$ in the notation of Lemma \ref{long first row},
which applies to settle this case.  So assume that
$X'$ has strictly more than $m$ columns.  Let $s'$ be a corner cell
of $X'$ of the form $(r,c)$ where $c>m$.  Note that $s'$
is also a corner cell of $X$, $S'$ and $S$.  Let
$Y'\ge X'$, $Y\ge X$, $T'\ge S'$ and $T\ge S$ be the cyclages
that end at $s'$.  These cyclages exist, since
the bumping paths of the reverse row insertions
on the tableaux $X',X,S',S$ starting at $s'$, move weakly east
so that the cycled letters are always strictly greater than one.
It is not hard to see that the cocyclage on $X'$ (resp. $X$)
at $s'$, commutes with the row insertion of the word $1^m$
(resp. $Z_m$).  So $T' = P(Y' 1^m)$, $T=P(Y Z_m)$.
By induction, the lemma can be
applied to $Y',Y,T',T$ so that $\theta'_{\muhat}(Y')=Y$ if and only if
$\theta_\mu(T')=T$.  But the $\theta$ operators commute with cyclage,
so we are done.
\end{proof}

Finally the proof of Theorem \ref{cyc image} is given.
\begin{proof} In view of Propositions \ref{cat and type}
and \ref{row and col cat}, it is enough to show that
the standard tableau $S$ is in the image of $\theta_\mu$,
if and only if $S$ is $\mu$-column catabolizable.

By Lemma \ref{long first row} we may assume that
$S$ contains $Z_m$.  Let $S_e$ and $S_w$
(resp. $S'_e$ and $S'_w$) be the east
and west tableaux obtained by slicing the skew tableau
$S-Z_m$ (resp. $S'-(1^m)$) vertically between the $m$-th and
$(m+1)$-st columns.  Let $X'=P(S'_e S'_w)$ and $X=P(S_e S_w)$.
The conditions for Lemma \ref{theta and north hat} are satisfied,
so $\theta_\mu(S')=S$ if and only if
$\theta'_{\muhat}(X') = X$.  The following are equivalent.
\begin{enumerate}
\item $\theta_\mu(S')= S$.
\item $\theta'_{\muhat}(X') = X = \CCat_m(S)$.
\item $\CCat_m(S)$ is $\muhat$-column catabolizable.
\item $S$ is $\mu$-column catabolizable.
\end{enumerate}
$(1)\Leftrightarrow(2)$ has just been shown,
$(2)\Leftrightarrow(3)$ follows by induction, and
$(3)\Leftrightarrow(4)$ holds by definition.
\end{proof}

\subsection{The cocharge Kostka special case of Conjecture \ref{cat conj}}

\begin{prop} Conjecture \ref{cat conj} holds when
$\eta$ is a partition $\mu=(\mu_1,\mu_2,\dots,\mu_t)$
and $\gamma=(1^t)$.
\end{prop}
\begin{proof} Let $R=(R_1,\dots,R_t)$ be the sequence of
partitions corresponding to the pair of weights $\gamma$
and $\eta=\mu$.  Then $R_i=(1^{\mu_i})$, the partition
consisting of a single column of height $\mu_i$.
It is not hard to see that a standard tableau $S$ is
$R$-catabolizable if and only if the transpose tableau
$S^t$ is $\mu$-column catabolizable.  Also it is easy to
see that $\charge(S) = \cocharge(S^t)$.  The result follows from
the cocharge Kostka-Foulkes special case in subsection
\ref{special case section} and the results of the previous section.
\end{proof}

\end{document}